# DIRECTIONS AND PROJECTIVE SHAPES[1]

By Kanti V. Mardia and Vic Patrangenaru

*University of Leeds and Texas Tech University*

This paper deals with projective shape analysis, which is a study of finite configurations of points modulo projective transformations. The topic has various applications in machine vision. We introduce a convenient projective shape space, as well as an appropriate coordinate system for this shape space. For generic configurations of $k$ points in $m$ dimensions, the resulting projective shape space is identified as a product of $k - m - 2$ copies of axial spaces $\mathbb{R}P^m$. This identification leads to the need for developing multivariate directional and multivariate axial analysis and we propose parametric models, as well as nonparametric methods, for these areas. In particular, we investigate the Frećhet extrinsic mean for the multivariate axial case. Asymptotic distributions of the appropriate parametric and nonparametric tests are derived. We illustrate our methodology with examples from machine vision.

**1. Introduction.** Consider a configuration of points in $\mathbb{R}^m$. "Shape" deals with the residual structure of this configuration when certain transformations are filtered out. More specifically, the shape of a configuration consists of its equivalence class under a group of transformations. Important groups for machine vision are the similarity group, the affine group and the projective group. Here the group action describes the way in which an image is captured. For instance, if two different images of the same scene are obtained using a pinhole camera, the corresponding transformation between the two images is the composition of two central projections, which is a projective transformation. If the two central projections can be approximated by parallel projections, which is the case of remote views of the same planar scene, the projective transformation can be approximated by an affine transformation. Further, if these parallel projections are orthogonal projections on the

Received October 2002; revised November 2004.
[1]Supported in part by NSA Grant 904-02-1-0082 and NSF Grant DMS-04-06151.

plane of the camera, this affine transformation can be approximated by a similarity transformation. Therefore, the relationships between these shapes are as follows: if two configurations have the same similarity shape, then they automatically have the same affine shape; if they have the same affine shape, they will have the same projective shape. For example, two squares of different sizes have the same similarity, affine and projective shape, whereas a square and a rectangle have the same affine and projective shape but not the same similarity shape. On the other hand, a square and a kite have the same projective shape but not the same affine shape.

The word "shape" often refers in statistics to similarity shape where only the effects of translation, scale and rotation have been filtered out (see, e.g., [4, 5]). In recent years, substantial progress has been made in similarity shape analysis, since appropriate shape space (e.g., Kendall's space) and shape coordinates (e.g., Bookstein coordinates) have been available. A simple example of Bookstein coordinates is for the shape of a triangle where the shape coordinates are obtained after taking one of the vertices as the origin and rotating the triangle so that the base of the triangle lies on the $x$-axis, and then rescaling the base to the unit size. The motivation behind such coordinate systems is similar to those in directional statistics, where to analyze spherical data one requires a coordinate system such as longitude and latitude (see, e.g., [18]).

Affine shape has also received some attention; see, for example, [10, 24]. Sparr [24] has shown that the space of affine shapes is a Grassmann manifold. For affine shape in 2-D, we can obtain shape coordinates by using three points that determine the direction and the origin of the axes, and the unit length between the points on each of these two axes.

Progress in projective shape analysis has been somewhat slow by not having a convenient shape space, though considerable work has appeared on projective invariants (see, e.g., [11, 12, 20]). We propose a convenient projective shape space, as well as an appropriate coordinate system for this shape space.

The plan of the paper is as follows. In Section 2 we propose our approach in "projective shape analysis," which has its basis on the idea of constructing a *projective frame* selected from the points of a generic configuration. The resulting projective shape space is a product of $k - m - 2$ copies of axial spaces $\mathbb{R}P^m$. This axial representation leads to various questions in multivariate directional statistics. To address these questions, in Section 3 we first discuss some parametric models, especially for the multivariate circular case. As a starting point, we consider certain von Mises circular distributions. These provide good approximations to marginal distributions on the circle of cross-ratios with normal errors at landmarks (i.e., offset projective distributions), as argued through simulations in [11]. We then treat the case of concentrated data by using a directional representation in a tangent space. In particular, the procedure is illustrated by constructing a two-sample test.

In Section 4 we consider estimation of certain means, both asymptotically and through bootstrap methods. In particular, we treat the multivariate axial case, highlighting the extrinsic mean; for $m = 1$, the circular extrinsic mean is well studied (see, e.g., [18]) and is generally referred to as mean direction. Theorem 4.1 provides asymptotic distributions of certain test statistics required for the estimation of the extrinsic mean of projective shapes for any $m$, and Corollary 4.1 provides some bootstrap approximations for these asymptotic distributions. In Section 4 we also provide a two-sample test for extrinsic means of projective shape through a bootstrapping result.

In Section 5 we illustrate our methodology through three examples in object recognition. The first two examples concern building recognition and use circular and univariate spherical statistics, whereas the third example is about face recognition and uses bivariate spherical statistics. Of course, the realm of applications is much wider, covering other types of multivariate axial data. In Section 6 we present a strategy for general statistical shape analysis where the shapes are regarded as orbits of certain Lie group actions on a direct product of a number of copies of a manifold.

**2. The projective shape space.** Recall that the real projective space in $m$ dimensions, $\mathbb{R}P^m$, is the set of axes going through the origin of $\mathbb{R}^{m+1}$. If $X = (X^1, \ldots, X^{m+1}) \in \mathbb{R}^{m+1} \setminus \{0\}$, then

$$[X] = [X^1 : X^2 : \ldots : X^{m+1}] = \{\lambda X, \lambda \neq 0\}$$

is a *projective point* in $\mathbb{R}P^m$; we will reserve the notation $[\cdot]$ for the projective points throughout. In an alternative description, a point $p \in \mathbb{R}P^m$ is given by $p = [z^1 : z^2 : \ldots : z^{m+1}]$, where

$$(z^1)^2 + (z^2)^2 + \cdots + (z^{m+1})^2 = 1.$$

A linear variety $v$ of dimension $d$ is given by $v = \{[x], x \in V \setminus 0\}$, where $V$ is a $(d+1)$-dimensional vector subspace of $\mathbb{R}^{m+1}$. In particular, a projective line $l$ is a set associated with a vector plane $V$ in $\mathbb{R}^{m+1}$, $l = \{[x], x \in V \setminus 0\}$. A number of points in $\mathbb{R}P^m$ are collinear if they lie on a projective line.

The Euclidean space $\mathbb{R}^m$ can be embedded in $\mathbb{R}P^m$, preserving collinearity. Such a standard *affine* embedding, missing only a hyperplane at infinity, is

$$x = (x^1, \ldots, x^m) \to [x^1 : \ldots : x^m : 1].$$

This leads to the notion of *affine* or *inhomogeneous* coordinates of a point

$$p = [X] = [X^1 : \ldots : X^m : X^{m+1}], \qquad X^{m+1} \neq 0,$$

to be defined as

$$(x^1, x^2, \ldots, x^m) = \left(\frac{X^1}{X^{m+1}}, \ldots, \frac{X^m}{X^{m+1}}\right),$$

as opposed to the *homogeneous* coordinates of $p$, $(X^1, \ldots, X^{m+1})$, which are defined up to a multiplicative constant only. However, the coordinates of interest in projective shape analysis are neither affine nor homogeneous. We need coordinates that are invariant with respect to the group of projective (general linear) transformations $PGL(m)$. A *projective transformation* $\alpha$ of $\mathbb{R}P^m$ is defined in terms of an $(m+1) \times (m+1)$ nonsingular matrix $A \in GL(m+1, \mathbb{R})$ by

$$\alpha([X^1 : \ldots : X^{m+1}]) = [A(X^1, \ldots, X^{m+1})^T].$$

The linear span of a subset of $\mathbb{R}P^m$ is the smallest linear variety containing that subset. Note that $k$ points in $\mathbb{R}P^m$ with $k \geq m+2$ are in *general position* if their linear span is $\mathbb{R}P^m$.

DEFINITION 2.1. A *projective frame* in $\mathbb{R}P^m$ is an ordered system of $m+2$ points in general position.

In computer vision, a projective frame is called a *projective basis* by some authors (Heyden [13], page 8; Faugeras and Luong [6], page 81). Let $(e_1, \ldots, e_{m+1})$ be the standard basis of $\mathbb{R}^{m+1}$. The *standard* projective frame is $([e_1], \ldots, [e_{m+1}], [e_1 + \cdots + e_{m+1}])$. The last point of this frame is referred to as the *unit* point.

PROPOSITION 2.1. *Given two projective frames $\pi_1 = (p_{1,1}, \ldots, p_{1,m+2})$ and $\pi_2 = (p_{2,1}, \ldots, p_{2,m+2})$, there is a unique $\beta \in PGL(m)$ with $\beta(p_{1,j}) = p_{2,j}$, $j = 1, 2, \ldots, m+2$.*

A proof follows on noting that, given a projective frame $\pi = (p_1, \ldots, p_{m+2})$, there is a unique $\alpha \in PGL(m)$ with

(2.1)  $\alpha([e_j]) = p_j, \qquad j = 1, \ldots, m+1, \qquad \alpha([e_1 + \cdots + e_{m+1}]) = p_{m+2}.$

REMARK 2.1. If $k > m+2$, we consider the set $\mathcal{FC}_m^k$ consisting of configurations of points $(p_1, \ldots, p_k)$ for which there is a subset of indices $i_1 < \cdots < i_{m+2}$ such that $(p_{i_1}, \ldots, p_{i_{m+2}})$ is a projective frame.

From Proposition 2.1, $\mathcal{FC}_m^k$ is an invariant generic subset of $\mathcal{C}_m^k$. It can be shown by considering, for example, $m = 1, k = 4$, that the corresponding shape space is a manifold. Let us denote the projective shape of $(p_1, p_2, p_3, p_4) \in \mathcal{FC}_1^4$ by $\sigma$.

Any projective shape is in one of the sets $U_{123}, U_{124}, U_{134}$ or $U_{234}$, where, for $i < j < k$,

(2.2) $\qquad U_{ijk} = \{\sigma | (p_i, p_j, p_k) \text{ is a projective frame}\}.$

Assume $p_r = [x_r : 1], r = 1, 2, 3, 4$. Then, from (2.2) we may define the charts $\psi_{ijk}$ by $\psi_{ijk}(\sigma) = c(p_i, p_j, p_k, p_l)$, where $c(\cdot)$ is a cross-ratio defined by

$$c(p_i, p_j, p_k, p_l) = \{(x_i - x_k)(x_l - x_j)\} / \{(x_i - x_j)(x_l - x_k)\}$$

and $\{i, j, k, l\} = \{1, 2, 3, 4\}$. On permuting indices, we find that

$$\psi_{124} = \frac{1}{\psi_{123}}, \qquad \psi_{134} = 1 - \psi_{124}, \qquad \psi_{234} = \frac{\psi_{134}}{\psi_{134} - 1}.$$

Since the transition maps between these charts are differentiable, we conclude that the projective shape space associated with $\mathcal{FC}_1^4$ is a 1-D manifold. Projective shape spaces associated with $\mathcal{FC}_m^k$ are expected to be more complicated. However, we prefer to restrict ourselves to a subset of generic configurations, such that the corresponding shape space has a natural structure of symmetric space (see Proposition 2.3), and for which the computations can be carried out using standard statistical packages.

DEFINITION 2.2. The *axis* of a point $p \in \mathbb{R}P^m$ with respect to a projective frame $\pi = (p_1, \ldots, p_{m+2})$ is defined as $p^\pi = \alpha^{-1}(p)$, where $\alpha \in PGL(m)$ is given by (2.1). A geometric interpretation of $p^\pi$ is given below.

Assume $f_{m+2} \in \mathbb{R}^{m+1}$ is a representative of $p_{m+2}$. Since $(p_1, \ldots, p_{m+2})$ are in general position, $f_{m+2}$ can be written in a unique way as a sum $f_{m+2} = f_1 + \cdots + f_{m+1}$, where $[f_j] = p_j$, for $j = 1, \ldots, m+1$. The vectors $f_1, \ldots, f_{m+1}$ form a basis of $\mathbb{R}^{m+1}$, and let $f \in \mathbb{R}^{m+1}$ be a representative of $p$. Then we denote $Y^1, \ldots, Y^{m+1}$ as the components of $f$ with respect to this basis. Note that since the selection of $f_{m+2}$ and of $f$ is unique up to a multiplicative constant, the projective point $[Y^1 : \ldots : Y^{m+1}]$ is well defined in $\mathbb{R}P^m$ and $p^\pi = [Y^1 : \ldots : Y^{m+1}]$. This representation of projective coordinates is displayed in Figure 1 for $m = 1$. Figure 1(a) constructs the coordinates $Y_1, Y_2$ of $f$ with respect to the frame $([f_1], [f_2], [f_3])$; Figure 1(b) shows the corresponding projective point $[Y_1 : Y_2]$.

Let us assume that $x_1, \ldots, x_{m+2}$ are points in general position and let $x = (x^1, \ldots, x^m)^T$ be an arbitrary point in $\mathbb{R}^m$. In this notation the axis of $x$ with respect to the projective frame associated with $m+2$ points $x_1, \ldots, x_{m+2}$ is the same as the axis of $p = [x^1 : \ldots : x^m : 1]$ with respect to $(p_1, \ldots, p_{m+2})$. Using the above geometric interpretation, we determine the axis of $x$ in the following proposition.

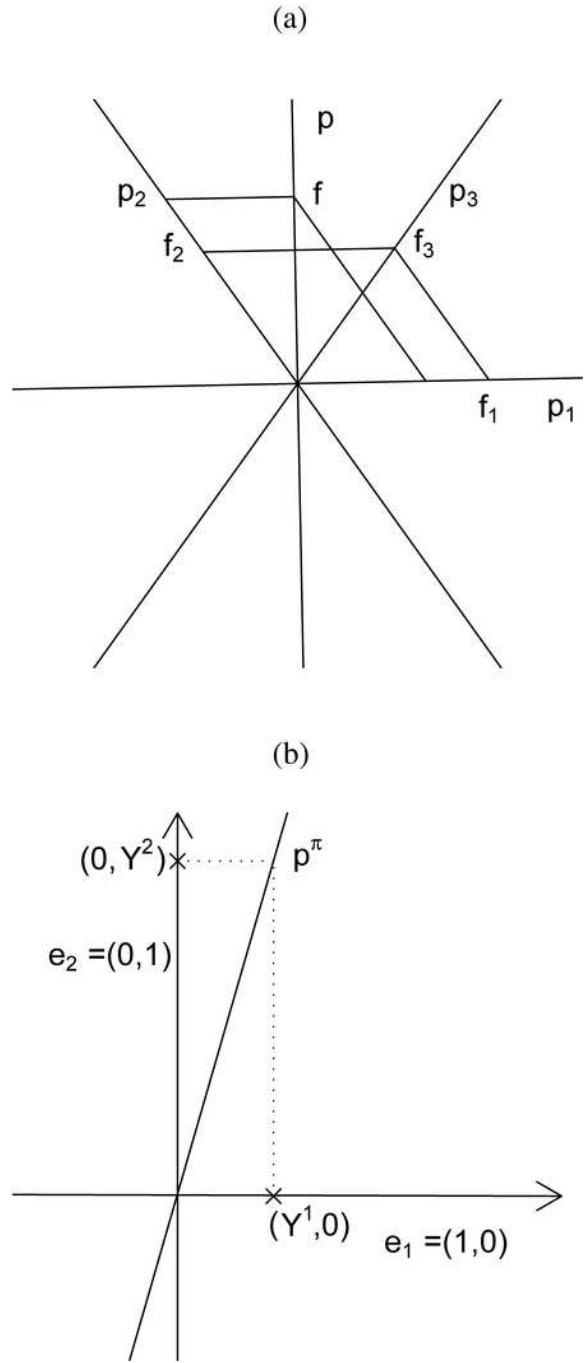

FIG. 1. *Projective coordinates for $m = 1$. (a) Projective frame $\pi = (p_1, p_2, p_3)$ and a projective point p, (b) projective coordinates of p with respect to $\pi$.*

PROPOSITION 2.2. *The projective axis of a point $x$ is given by*

(2.3) $$[z^1(x):z^2(x):\ldots:z^{m+1}(x)],$$

*where*

(2.4) $$z^j(x) = y^j(x)/\|y(x)\|, \qquad j=1,\ldots,m+1,$$

$$y^j(x) = v^j(x)/v^j(x_{m+2}),$$

(2.5) $$y(x)^T = (y^1(x),\ldots,y^{m+1}(x)),$$

$$z(x) = (z^1(x),\ldots,z^{m+1}(x))^T, \qquad \|z(x)\|^2 = 1$$

*and*

(2.6) $$v(x) = (v^1(x),\ldots,v^{m+1}(x))^T = U_m^{-1}p(x),$$

*with the $(m+1) \times (m+1)$ matrix*

(2.7) $$U_m = [p(x_1),\ldots,p(x_{m+1})],$$

*where $p(x) = (x^T, 1)^T$.*

PROOF. We note that there is a unique vector $\beta \in \mathbb{R}^{m+1}$, $\beta^T = (\beta_1,\ldots,\beta_{m+1})$, such that $\beta_1 p(x_1) + \cdots + \beta_{m+1} p(x_{m+1}) = p(x_{m+2})$. Namely, $U_m \beta = p(x_{m+2})$ or $\beta = U_m^{-1} p(x_{m+2}) \equiv v(x_{m+2})$, so that

(2.8) $$\beta_i = v^i(x_{m+2}), \qquad i = 1,\ldots,m+1.$$

Let $A$ be the matrix

(2.9) $$A = U_m \operatorname{diag}(\beta).$$

If $(e_1,\ldots,e_{m+1})$ is the standard basis for $\mathbb{R}^{m+1}$, then $Ae_i = \beta_i p(x_i)$, $i = 1,\ldots,m+1$, and $A(e_1 + \cdots + e_{m+1}) = p(x_{m+2})$. This means that the projective transformation $\alpha$, given by $\alpha([x]) = [Ax]$, has the properties that $\alpha[e_i] = [p(x_i)], i = 1,\ldots,m+1$, and $\alpha[e_1 + \cdots + e_{m+1}] = [p(x_{m+2})]$. Hence, from Definition 2.2 it now follows that the homogeneous projective coordinates of $[p(x)]$ are given by

(2.10) $$y(x) = A^{-1}p(x).$$

From (2.8) and (2.9) we have

$$y(x) = \operatorname{diag}\left(\frac{1}{\beta_1},\ldots,\frac{1}{\beta_{m+1}}\right) U_m^{-1} p(x)$$

$$= \operatorname{diag}\left(\frac{1}{v^1(x_{m+2})},\ldots,\frac{1}{v^{m+1}(x_{m+2})}\right) U_m^{-1} p(x).$$

Hence, using the definition of $v(x)$ given by (2.6), (2.5) follows. □

REMARK 2.2. We will say that $y(x)$ are the projective coordinates of $x$ with respect to the projective frame generated by $(x_1, \ldots, x_{m+2})$, and note that $[z(x)]$ defined by (2.3) is the corresponding point on $\mathbb{R}P^m$.

REMARK 2.3. Note that we have $v(x_i) = e_i, i = 1, \ldots, m+1$.

Let $G(k, m)$ denote the set of all ordered systems of $k$ points $(p_1, \ldots, p_k)$ for which $(p_1, \ldots, p_{m+2})$ is a projective frame, $k > m+2$. $PGL(m)$ acts on $G(k, m)$ by $\alpha(p_1, \ldots, p_k) = (\alpha p_1, \ldots, \alpha p_k)$.

DEFINITION 2.3. The projective shape space $P\Sigma_m^k$ or space of *projective k–ads* in $\mathbb{R}P^m$ is the quotient $G(k,m)/PGL(m)$.

PROPOSITION 2.3. $P\Sigma_m^k$ *is a manifold diffeomorphic with* $(\mathbb{R}P^m)^q$, *where* $q = k - m - 2$.

PROOF. We define $F: P\Sigma_m^k \to (\mathbb{R}P^m)^q$ by

$$F((p_1, \ldots, p_k) \mod PGL(m)) = (p_{m+3}^\pi, \ldots, p_k^\pi),$$
(2.11)
$$\pi = (p_1, \ldots, p_{m+2}), p_i^\pi = [z^1(x_i) : \ldots : z^{m+1}(x_i)],$$

where $z(x_i) = (z^1(x_i), \ldots, z^{m+1}(x_i))^T, \|z(x_i)\| = 1, i = m+3, \ldots, k$, and $z(\cdot)$ is given by (2.3).

The mapping $F$ is a well-defined diffeomorphism between $P\Sigma_m^k$ and a product of real projective spaces.

Note that (2.11) defines an *axial* representation of the projective shape. In this representation, for $m = 1$ we can write (2.11) as $p_j^\pi = [e^{i\phi_j}]$, where $\phi_j$ is the angle of an axis. Then doubling $\phi_j$ takes us to an oriented direction $e^{i\theta_j} \in S^1$. Further, we assume that $x_1, x_2, x_3$ yield a projective frame $\pi$ and $[x:1]$ is an arbitrary point on the projective line. Following the above algorithm for projective coordinates, from (2.5) and (2.6) we get

$$v^1(x) = \frac{x - x_2}{x_1 - x_2}, \qquad v^2(x) = \frac{x_1 - x}{x_1 - x_2},$$
$$y^1(x) = v^1(x)/v^1(x_3), \qquad y^2(x) = v^2(x)/v^2(x_3).$$

Thus, from (2.4) we have

$$z^1(x) = \frac{y^1(x)}{\{(y^1(x))^2 + (y^2(x))^2\}^{1/2}}, \qquad z^2(x) = \frac{y^2(x)}{\{(y^1(x))^2 + (y^2(x))^2\}^{1/2}}$$

or, equivalently,

(2.12) $$z^1(x) = \frac{x - x_2}{(x_3 - x_2)\, d(x)}, \qquad z^2(x) = \frac{x_1 - x}{(x_1 - x_3)\, d(x)},$$

where
$$d(x) = \left\{ \left( \frac{x - x_2}{x_3 - x_2} \right)^2 + \left( \frac{x_1 - x}{x_1 - x_3} \right)^2 \right\}^{1/2}.$$

That is, we can write

(2.13) $\quad p^\pi = [z^1(x) : z^2(x)], \qquad z^1(x) = \cos \phi(x), z^2(x) = \sin \phi(x).$

This representation of projective coordinates is displayed in Figure 2. Note that on eliminating $x$ from $y^1(x)$ and $y^2(x)$, we get $(x_3 - x_2)y^1(x) + (x_1 - x_3)y^2(x) = x_1 - x_2$. Since $x_1 \neq x_2 \neq x_3$, this equation of a line in the plane $(y^1(x), y^2(x))$ confirms that the angle $\phi(x)$ lies between 0 and $\pi$. $\square$

EXAMPLE 2.1. We now illustrate our coordinate system for a problem studied in machine vision by Heyden [13], pages 33–34. He has considered a configuration of five points from two images of a rectangular sheet of paper. The five points form a cross. For the first image the coordinates are

$$x_1 = (69, 53), \qquad x_2 = (591, 33), \qquad x_3 = (626, 402),$$
$$x_4 = (69, 430), \qquad x_5 = (344, 322).$$

The first four points are the corners and the last one is the center of the cross. Here $m = 2$ and, from (2.7) we have, on registering with respect to the frame $x_1$, $x_2$ and $x_3$:

$$U_2 = \begin{bmatrix} 69 & 591 & 626 \\ 53 & 33 & 402 \\ 1 & 1 & 1 \end{bmatrix}.$$

Hence, from (2.6) we find that $v(x_4)^T = (1.0683, -1.0862, 1.0180)$ and $v(x_5)^T = (0.5057, 0.0095, 0.4848)$. Thus, from (2.4) and (2.5) we get the following spherical representation of the projective shape:

$$z(x_5) = [0.7050 : -0.0131 : 0.7092] = z_1, \qquad \text{say.}$$

Similarly, for the second image we have $x_1 = (334, 69)$, $x_2 = (732, 290)$, $x_3 = (428, 504)$, $x_4 = (43, 200)$, $x_5 = (373, 243)$, leading to

$$z(x_5) = [0.7074 : -0.0060 : 0.7067] = z_2, \qquad \text{say.}$$

We will return to these coordinates in Example 4.1.

We now describe the alternative representation of projective shape due to Goodall and Mardia [11] and show the connection between the two representations. In their representation the projective shape of $(p_1, \ldots, p_k) \in G(k, m)$ is uniquely determined by its projective invariants. In fact, the projective coordinates, with respect to $(p_1, \ldots, p_{m+2})$ and the projective invariants $(\iota_{ji})$,

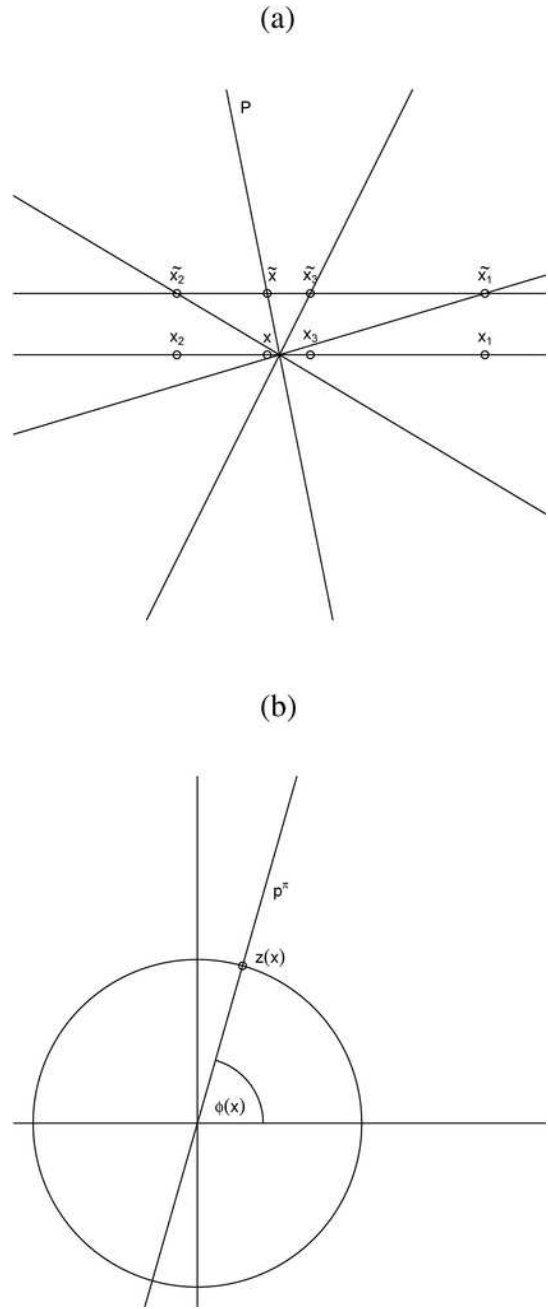

FIG. 2. *Projective coordinate of a point on a Euclidean line* (a) *shows the embedding of an affine line in the projective line and the projective points corresponding to the points on the real line;* (b) *shows the projective coordinate of the point $p = [x:1] = [\tilde{x}]$ with respect to the projective frame $p_1 = [\tilde{x}_1], p_2 = [\tilde{x}_2], p_3 = [\tilde{x}_3]$ and the corresponding angle of this projective coordinate, $\phi(x)$, which should be doubled to get the point $\theta(x)$ on the circle.*

$j = 1, \ldots, m; i = m+3, \ldots, k$, determine each other by the *invariant* representation

$$(2.14) \qquad p_i^\pi = [\iota_{1i} : \iota_{2i} : \ldots : \iota_{mi} : 1],$$

so that the $i$th projective coordinate is $\iota_{\cdot, i}$. When $m = 1, 2$, let us consider their invariant representation. For $m = 1$, from (2.12) and (2.13) we have

$$(2.15) \qquad p^\pi = [x:1]^\pi = \left[\frac{x-x_2}{x_3-x_2} : \frac{x_1-x}{x_1-x_3}\right] = [c(x_1, x_2, x_3, x) : 1].$$

This equation shows that the projective coordinate of $[x:1]$ with respect to $\pi$ in our representation [viz., equation (2.13)] and the cross-ratio $c(x_1, x_2, x_3, x)$ determine each other. Here $x$ could be any of the $k-3$ points. For $m = 2$ we now assume that $x_1, x_2, \ldots, x_k$ are $k$ points in $\mathbb{R}P^2, k > 4$. Let $z_i^j = z^j(x_i), i = 5, \ldots, k$. Then the invariants $\iota_{1i}$ and $\iota_{2i}$ in (2.14) are the cross-ratios

$$(2.16) \qquad \iota_{1i} = \frac{z_i^1}{z_i^3}, \qquad \iota_{2i} = \frac{z_i^2}{z_i^3},$$

where $\iota_{1i}$ and $\iota_{2i}$ are the cross-ratios determined on a transversal by the pencil of lines joining the point $x_1$ to the points $x_2, x_3, x_4, x_i$ and by the pencil of lines joining the point $x_2$ to the points $x_1, x_3, x_4, x_i$, respectively. There are some parallel ideas between coordinates for similarity shape, affine shape and the projective shape considered here. We now give registration frames for the three shapes in 2-D. Let us consider planar similarity shape. The shape can be registered using two points; for example, we can use the registration frame with the points $(0,0)$ and $(1,0)$. Figure 3(a) shows the original configuration for $k = 4$ and Figure 3(b) shows its Bookstein coordinates. Figure 3(c) is discussed below. In the case of affine shape ("intermediate" between similarity shape and projective shape), we can choose the registration frame consisting of the three points $(0,0), (1,0)$ and $(0,1)$. In projective shape in 2-D, using homogeneous coordinates, we can select the registration frame consisting of the points $(0,0), (1,0), (0,1)$ and $(1,1)$. In inhomogeneous coordinates, the registration frame corresponds to the points $(1,0,0), (0,1,0), (0,0,1)$ and $(1,1,1)$ in 3-D. The steps in projective shape registration for 1-D and 2-D are as follows:

*Case* $m = 1$.

1.0. Start with a configuration of $k$ points.

1.1. Register each of the last $k-3$ points with respect to first three points leading to (2.6).

1.2. Transform these registered points by (2.4), leading to one point on the Cartesian product $\mathbb{R}P^1 \times \cdots \times \mathbb{R}P^1$ of $k-3$ copies of $\mathbb{R}P^1$, that is, the projective shape of a linear configuration of $k$ points is equivalent to $k-3$ axes in $\mathbb{R}^2$.

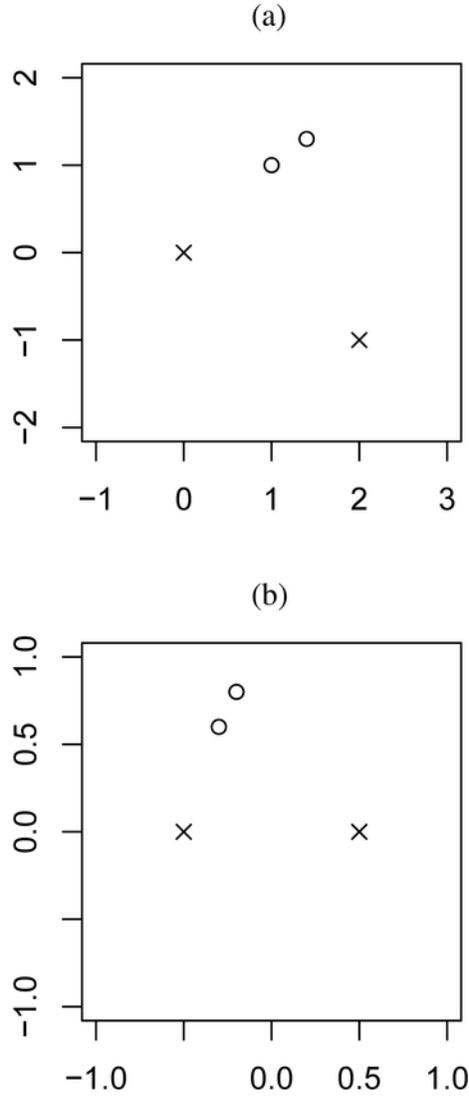

Fig. 3. (a) *Similarity shape: four points in* 2-*D with the two base points marked by* "×," *and the other by* "∘." (b) *Bookstein registration for Figure* 3(a) *with respect to the frame* $(-\frac{1}{2}, 0)$ *and* $(\frac{1}{2}, 0)$.

1.3. Transform these $k-3$ axes to directions by doubling the angles. Thus, we get an observation on a $(k-3)$-dimensional torus.

Figure 3(c) gives a schematic diagram of these first two steps for $k = 5$.

*Case* $m = 2$.

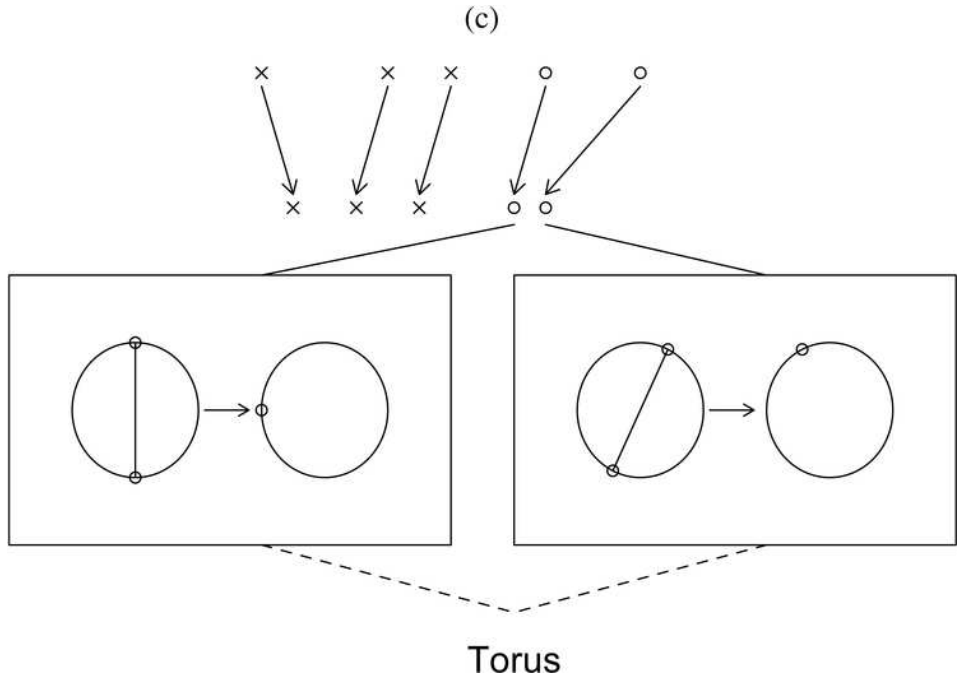

Fig. 3. (*Continued*). (c) *Projective shape: schematic diagram for five points in* $\mathbb{R}^1$ *and their registration steps.*

2.0. Start with a configuration of $k$ points.

2.1. Register the projective coordinates of the last $k-4$ points with respect to the first four points, say.

2.2. Transform these registered points leading to one point on the Cartesian product $\mathbb{R}P^2 \times \cdots \times \mathbb{R}P^2$ of $k-4$ copies of $\mathbb{R}P^2$, that is, the projective shape of a planar configuration of $k$ points is equivalent to $k-4$ axes in $\mathbb{R}^3$.

These schematic constructions are displayed in Figure 4 for $m=2$ and $k=6$. Note that unlike the case $m=1$, for $m=2$, $\mathbb{R}P^2$ cannot be visualized in three dimensions, as we need at least four dimensions to immerse the real projective plane into a Euclidean space without double points (see [25], pages I-9–I-13). A rigorous geometric construction for $m=1$ and $k=4$ has already been given in Figures 1 and 2.

In the subsequent discussion, we will work on the axial $z$-coordinates derived from the projective coordinates $x$ with respect to the frame $(x_1, \ldots, x_{m+2})$. From now on the $z$-coordinates will be written as $[X] = ([x_1], \ldots, [x_q])$, which corresponds to an $(m+1) \times q$ matrix. In the case of the directional representation (for concentrated data), we will write $X = (x_1, \ldots, x_q)$ when there is no ambiguity.

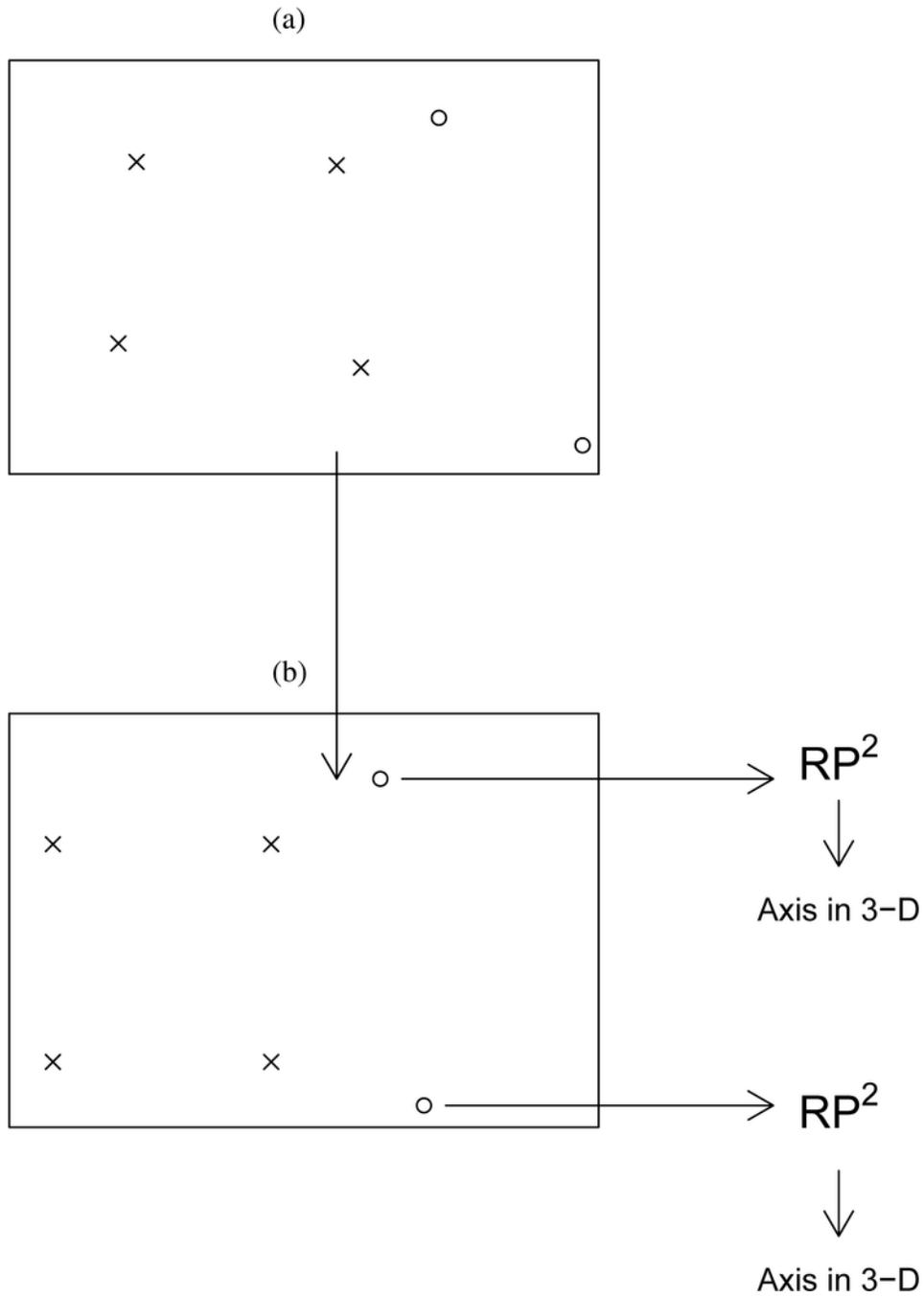

FIG. 4. (a) *Six points in* $\mathbb{R}^2$ *and* (b) *a schematic representation of their projective coordinates.*

**3. Spherical distributions for projective shape.** We have shown in Proposition 2.3 that the projective shape of $k$ landmarks in $m$ dimensions can be represented by $q$ axial variables given by (2.3). Thus, we could use appropriate axial distributions developed in directional statistics as models for projective shape (see, e.g., [18]). Simulation studies performed by Goodall and Mardia [11] suggest that von Mises distributions are appropriate approximations to the angular version of the cross-ratios under isotropic normal variation at landmarks; see their chi-square plots in Figures 9 and 10. Also, various applications are given in that paper. That is, for $m = 1$ and $k = 4$, the angle $\theta = 2\phi$ can be considered to have the von Mises distribution with probability density function

$$f(\theta; \mu, \kappa) = \{2\pi I_0(\kappa)\}^{-1} \exp\{\kappa \cos(\theta - \mu)\}, \qquad \theta \in (0, 2\pi], \kappa \geq 0,$$

where $\mu \in (0, 2\pi]$ is the mean direction, $\kappa$ is the concentration parameter and $I_0(\kappa)$ is the modified Bessel function of the first kind and order zero.

In general, for $m = 1$ and any number of landmarks $k$, we require a multivariate von Mises distribution on the torus (after doubling the angles). Mardia [17] proposed the following family of bivariate von Mises densities:

$$\text{const.} \times \exp\{\kappa_1 \cos\theta_1 + \kappa_2 \cos\theta_2 + \kappa_3 \cos\theta_1 \cos\theta_2$$
$$+ \kappa_4 \cos\theta_1 \sin\theta_2 + \kappa_5 \sin\theta_1 \cos\theta_2 + \kappa_6 \sin\theta_1 \sin\theta_2\}.$$

The distribution could be generalized by mapping $\theta_1 \to \theta_1 - \mu$ and $\theta_2 \to \theta_2 - \nu$. A particular density of interest is proportional to

$$\exp\{\kappa_1 \cos(\theta_1 - \mu) + \kappa_2 \cos(\theta_2 - \nu) - \kappa_3 \cos(\theta_1 - \mu - \theta_2 + \nu)\},$$

where $-\pi < \theta_1, \theta_2 \leq \pi, -\pi < \mu, \nu \leq \pi, \kappa_1 \geq 0, \kappa_2 \geq 0$ and $\kappa_1 \geq \kappa_3 \geq 0, \kappa_2 \geq \kappa_3 \geq 0$. The density has the advantage that the parameters here have no redundancy for large concentration since the distribution tends to a full bivariate normal distribution. For another member of this family, see [23]. A multivariate extension of the distribution for variables $\theta_1, \ldots, \theta_q$ can be written down with density proportional to

$$(3.1) \quad \exp\{\Sigma a_s \cos\theta_s + \Sigma b_s \sin\theta_s + \Sigma a_{st} \cos\theta_s \cos\theta_t$$
$$+ \Sigma b_{st} \cos\theta_s \sin\theta_t + \Sigma c_{st} \sin\theta_s \sin\theta_t\},$$

where $a_{ss} = b_{ss} = c_{ss} = 0, b_{st} \neq b_{ts}$, and $s, t = 1, \ldots, q$.

For concentrated projective shape data of configurations of $k$ points in general position $(x_1, \ldots, x_k)$ where $(x_1, \ldots, x_{m+2})$ yields a projective frame in $\mathbb{R}^m$, one may simply consider a multivariate directional representation $(z(x_{m+3}), \ldots, z(x_k))$, where $z(x_i) \in S^m$ are defined by (2.3). If there is no ambiguity, we will denote the vector $(z(x_{m+3}), \ldots, z(x_k))$ as $(x_1, \ldots, x_q)$. A projective shape is represented as a point in a direct product of $q$ copies of $S^m$. We now consider models on $(S^m)^q$.

*Distributions in the tangent space.* For $m = 1$, let $\gamma_s$ be the mean direction of $\theta_s$, $s = 1, \ldots, q$. Then for $\gamma_s = 0, a_{st} = b_{st} = 0$, $b_s = 0$ and large $a_s$, it can be shown that $(\sin \theta_1, \ldots, \sin \theta_q)$ is $N(0, \Sigma)$, where $\Sigma$ is the "asymptotic" covariance which depends on the population parameters. If the angle $\gamma_s$ is not zero, then we have approximately

$$(3.2) \qquad (\sin(\theta_1 - \gamma_1), \ldots, \sin(\theta_q - \gamma_q))^T \sim N(0, \Sigma).$$

This approximation can be extended to spherical variables for $m > 1$. Let $X = (x_1, \ldots, x_q)$ be an $(m+1) \times q$ random matrix with $x_s \in S^m$. Note that these $x$'s are not to be confused with the notation elsewhere in the paper. Suppose that $\mu_s$ is the mean directional vector of $x_s, s = 1, \ldots, q$, and let $\mu = (\mu_1, \ldots, \mu_q)$ denote the population mean directional matrix of $X$. Then define the spherical tangent coordinates of $x_s$ by

$$(3.3) \qquad v_s = (I - \mu_s \mu_s^T) x_s, \qquad s = 1, \ldots, q.$$

For the circular case with $x_s^T = (\sin \theta_s, \cos \theta_s)$ and $\mu_s^T = (\sin \gamma_s, \cos \gamma_s)$, we have $v_s = \sin(\theta_s - \gamma_s) e_s$, where $e_s$ is a unit tangent vector at $\mu_s, s = 1, \ldots, q$. Then for concentrated data, it is plausible to extend (3.2) to the joint distribution of $v_1, \ldots, v_q$, which has a multivariate normal with zero mean vector and some covariance matrix. Such tangent approximations are commonly used in shape analysis (see, e.g., [5], Chapter 7) and found to be very effective in practice. We now give a few results without any proof since their proofs are similar to those in shape analysis.

Let $X_i, i = 1, \ldots, n$, be a random sample from the population with random matrix $X$, where $X_i = (x_{i,1}, \ldots, x_{i,q})$. We estimate the population multivariate mean directional matrix $\mu$ by the sample mean directional matrix $\hat{\mu} = (\hat{\mu}_1, \ldots, \hat{\mu}_q)$, where $\hat{\mu}_i, i = 1, \ldots, q$, are the standard spherical sample mean directions. Let $\tilde{v}_{i,s} = [I - \hat{\mu}_s \hat{\mu}_s^T] x_{i,s}, s = 1, \ldots, q$. Further, let $\bar{\tilde{v}}$ and $S$ be, respectively, the sample mean and covariance matrix of $\tilde{v}_i = (\tilde{v}_{i,1}^T, \ldots, \tilde{v}_{i,q}^T)^T, i = 1, \ldots, n$. Then the Mahalanobis distance squared $D^2 = \bar{\tilde{v}}^T S^- \bar{\tilde{v}}$ has rank $M = mq$, where $S^-$ is the Moore–Penrose generalized inverse of $S$.

PROPOSITION 3.1. *For concentrated data the approximate distribution of Hotelling's $T^2$ statistic is given by*

$$T^2 = \frac{n-M}{M} D^2 \sim F_{M, n-M},$$

*and the asymptotic distribution of $T^2$ is $\chi_M^2$, where $M = mq$.*

*Tangent space inference.* Using the above strategy, we now construct a two-sample Hotelling's $T^2$ test. Let $X_{1i}, i = 1, \ldots, n_1$, and $X_{2j}, j = 1, \ldots, n_2$, be two independent samples, where $X_{1i}$ and $X_{2j}$ are two $(m+1) \times q$ matrices where each column lies in $S^m$. Suppose that $\mu_1$ and $\mu_2$ are the respective mean population matrices. We wish to test the hypotheses

$$H_0: \mu_1 = \mu_2 \quad \text{vs.} \quad H_1: \mu_1 \neq \mu_2.$$

Let $\hat{\mu}$ be the matrix of the combined (sample) mean directions given by $(\hat{\mu}_1, \ldots, \hat{\mu}_q)$. Let

$$(3.4) \quad v_{i,s} = [I - \hat{\mu}_s \hat{\mu}_s^T] x_{1i,s}, \qquad w_{j,s} = [I - \hat{\mu}_s \hat{\mu}_s^T] x_{2j,s},$$

where $i = 1, \ldots, n_1, j = 1, \ldots, n_2, s = 1, \ldots, q$, and $X_{1i} = (x_{1i,1}, \ldots, x_{1i,q})$, $X_{2j} = (x_{2j,1}, \ldots, x_{2j,q})$. Assuming that these two independent samples are from the normal populations in this tangent space with the same covariance matrix, we find that the Mahalanobis distance squared between $\bar{v}$ and $\bar{w}$ is

$$D^2 = (\bar{v} - \bar{w})^T S^- (\bar{v} - \bar{w}),$$

where $\bar{v}$ is the $mq \times 1$ vector of the means of $(v_{i,1}^T, \ldots, v_{i,q}^T)^T, i = 1, \ldots, n_1$. The mean vector $\bar{w}$ is similarly defined for the second sample. Further, $S^-$ is the Moore–Penrose generalized inverse of $S = (n_1 S_1 + n_2 S_2)/(n_1 + n_2 - 2)$, where $S_1$ and $S_2$ are the sample covariance matrices. Note that the rank of $S$ is $M$.

PROPOSITION 3.2. *For concentrated data the two-sample Hotelling's $T^2$ statistic is given by*

$$(3.5) \quad T^2 = \frac{n_1 n_2 (n_1 + n_2 - M - 1) D^2}{(n_1 + n_2)(n_1 + n_2 - 2) M},$$

*and under $H_0$, $T^2$ is approximately distributed as $F_{M, n_1 + n_2 - M - 1}$. A general distribution for nonconcentrated axial data could be taken as a multivariate Bingham distribution.*

**4. An extrinsic mean.** Here, we define an appropriate location parameter for a probability distribution on $P\Sigma_m^k$. Then for a random sample $(X_1, \ldots, X_n)$ from a given probability distribution $Q$, we find a consistent estimator for the population location and derive its asymptotic distribution. We assume that a distance $\rho$ on $P\Sigma_m^k$ is specified. Through this distance $\rho$ it is possible to define an index of dispersion. Given a probability measure $Q$ on $P\Sigma_m^k$, following the general treatment of Fréchet [9] (see also [15]), we define, for $y \in P\Sigma_m^k$, the function

$$(4.1) \quad \mathcal{F}_Q(y) = E[\rho^2(X, y)] = \int_{P\Sigma_m^k} \rho^2(x, y) Q(dx).$$

Assume there is a unique $y \in P\Sigma_m^k$ such that $\inf_{\psi \in P\Sigma_m^k} \mathcal{F}_Q(\psi) = \mathcal{F}_Q(y)$; such a $y$ is said to be the *Fréchet population mean*, $y := \mu_\mathcal{F}$.

REMARK 4.1. For any probability measure $Q$ on $\mathbb{R}$ the mean is always unique. In general, this is true for a probability measure on $P\Sigma_m^k$, but there are some exceptions. For example, for the uniform distribution on $S^1 = P\Sigma_1^4$ with chord distance $\rho$ defined by

$$\rho^2(e^{i\theta}, e^{i\psi}) = 1 - \cos(\theta - \psi),$$

we have $F_Q(y) = 2\pi$ for any point $y$ on $S^1$; thus $y$ is not unique in this case.

Let $Y_1, Y_2, \ldots, Y_n$ be independent and identically distributed random variables with probability measure $Q$ and let $\widehat{Q}_n$ be the empirical probability measure

(4.2) $$\widehat{Q}_n = \frac{1}{n}(\delta_{Y_1} + \cdots + \delta_{Y_n}).$$

The *Fréchet sample mean set* is the set $\hat{\mu}_{\mathcal{F}} = \{\hat{y} \in P\Sigma_m^k, F_{\widehat{Q}_n}(\hat{y}) = \inf F_{\widehat{Q}_n}(y)\}$. If $\hat{\mu}_{\mathcal{F}}$ has a unique element, this element is called the Fréchet sample mean and is labelled $\bar{Y}_{\mathcal{F}}$. Further, Ziezold [29] has established the strong consistency of the Fréchet sample mean set on a compact metric space. Hence, from this result it follows that if $\mu_{\mathcal{F}}$ exists, then any measurable choice $\hat{y}$ from $\hat{\mu}_{\mathcal{F}}$ is a strongly consistent estimator of $\mu_{\mathcal{F}}$. If the mapping

$$j : P\Sigma_m^k \to E^N$$

is an embedding into a Euclidean space and we use the chord distance $\rho$,

$$\rho(p_1, p_2) = \|j(p_1) - j(p_2)\|, \qquad p_1, p_2 \in P\Sigma_m^k,$$

then the Fréchet mean is called the *extrinsic* mean and we denote it by $\mu_{E,j}$, or write simply $\mu_E$ when $j$ is known. One of the main features required in selecting an embedding is that the resulting mean is easily computable.

First, we consider the case of $q = 1$ in our formulation, which has already been studied in directional statistics (see [18, 21]). Note that for $m = 1$ the extrinsic mean of $Q$ for a population of projective shapes corresponds to the standard circular mean direction (see, e.g., [18], pages 29–30).

For $m > 1$, by Proposition 2.2 the space $P\Sigma_m^{m+3}$ is identified with the *axial space* $\mathbb{R}P^m$. We consider the embedding $j$ of $\mathbb{R}P^m$ into $S(m+1)$, the space of symmetric matrices [16] given for the directional representation $[x] = \{\pm x, \|x\| = 1\}$ by

(4.3) $$j([x]) = xx^T.$$

Here the Euclidean norm of a matrix $A \in S(m+1)$ is given by $\|A\|^2 = \operatorname{tr} AA^T$, that is, $A$ is an $(m+1) \times (m+1)$ symmetric matrix, and if $A = j([x])$ with $\|x\| = 1$, then $\|A\| = 1, A \geq 0$ and rank $A = 1$.

Let $[X], \|X\| = 1$ be a random vector in $\mathbb{R}P^m$. Then from [2] it follows that the extrinsic population mean exists if the largest eigenvalue of $E(XX^T)$ is simple (i.e., has multiplicity one). In this case $\mu_E = [\gamma]$, where $\gamma$ is an eigenvector of $E(XX^T)$ corresponding to the largest eigenvalue, with $\|\gamma\| = 1$. Moreover, if $[X_r]$, ($\|X_r\| = 1$), $r = 1, \ldots, n$, is a random sample from a probability measure $Q$ on $\mathbb{R}P^m$ and the extrinsic mean $\mu_E$ of $Q$ exists, then the extrinsic sample mean $\overline{[X]}_E$ is a strongly consistent estimator of $\mu_E(Q)$. Note that when it exists, $\overline{[X]}_E$ is given by

$$\overline{[X]}_E = [\mathbf{m}], \tag{4.4}$$

where $\mathbf{m}$ is a unit eigenvector of $\frac{1}{n}\sum_{r=1}^{n} X_r X_r^T$ corresponding to the largest eigenvalue. It may be noted that, in this case, $\overline{[X]}_E$ is also the maximum likelihood estimator (MLE) for the mean of a Bingham distribution ([16, 21]) and for the mean of the Dimroth–Watson distribution, whose density function at $[x]$ is proportional to $\exp(k(\mu \cdot x)^2)$, where $k$ is a constant. For these or more general parametric families, MLE asymptotics or bootstrap methods [7] are commonly used.

EXAMPLE 4.1. We find that the extrinsic sample mean shape of the two projective shapes given in Example 2.1 is given by

$$\overline{[z]}_E = [0.7062 : -0.0095 : 0.7080].$$

Since $n = 2$, we have $\overline{[z]}_E = [\bar{z}]_E = [(z_1 + z_2)/2]$. Heyden ([13], page 34) obtained the reconstructed point in our spherical coordinates as $[0.7070 : -0.0071 : 0.7071]$ by using a deterministic method. The two answers are very similar, though derived from two different methods.

The main argument for using the embedding $j$ in (4.3) is that, for this embedding, the extrinsic mean is easily computable via statistical packages with routines to carry out eigenanalysis. Another advantage is that the mapping $j$ is equivariant as shown in [16], and leads to the following multivariate extension to projective shapes.

If we use the axial representation of projective shapes, for $k \geq m + 3$ or $q \geq 1$, we can define an embedding $j_k$ of $(\mathbb{R}P^m)^q$ into $(S(m+1))^q$ in terms of $j$:

$$j_k([x_1], \ldots, [x_q]) = (j[x_1], \ldots, j[x_q]), \tag{4.5}$$

where $x_s \in \mathbb{R}^{m+1}, \|x_s\| = 1, s = 1, \ldots, q$. Again, it can be shown that if the largest eigenvalues of each of the $q$ matrix components of $E(j_k(Q))$ are simple, then the extrinsic mean $\mu_{j_k}(Q)$ exists and is given by

$$\mu := \mu_{j_k}(Q) = ([\gamma_1(m+1)], \ldots, [\gamma_q(m+1)]), \tag{4.6}$$

where $\gamma_s(m+1)$ is a unit eigenvector corresponding to the largest eigenvalue of the $s$th component of $E(j_k(Q))$. If $Y_r, r = 1, \ldots, n$, is a random sample from $Q$, then in the axial representation

$$(4.7) \qquad Y_r = ([X_{r,1}], \ldots, [X_{r,q}]), \qquad \|X_{r,s}\| = 1; s = 1, \ldots, q,$$

where $X_{r,s}$ is an $(m+1) \times 1$ vector and $Y_r$ is $(m+1) \times q$ matrix. Consider the matrix of sums of squares and products of entries of $X_{r,s}$ given by

$$(4.8) \qquad J_s = \frac{1}{n} \sum_{r=1}^{n} X_{r,s} X_{r,s}^T,$$

which is a well-defined $(m+1) \times (m+1)$ matrix. Let $d_s(a)$ and $g_s(a)$ be the eigenvalues in increasing order and the corresponding unit eigenvector of $J_s$, $a = 1, \ldots, m+1$. Then the extrinsic sample mean in this case is

$$(4.9) \qquad \bar{Y}_{n,E} = ([g_1(m+1)], \ldots, [g_q(m+1)]).$$

For $q = 1$, $\bar{Y}_{n,E}$ reduces to the mean given in (4.4), namely, $g_1(m+1) = \mathbf{m}$. Arrange the pairs of indices $(s, a), s = 1, \ldots, q, a = 1, \ldots, m$, in their lexicographic order, and define the $M \times M$ symmetric matrix $G$ by

$$(4.10) \quad \begin{aligned} &G_{(s,a),(t,b)} \\ &= n^{-1}(d_s(m+1) - d_s(a))^{-1}(d_t(m+1) - d_t(b))^{-1} \\ &\quad \times \sum_{r=1}^{n} (g_s^T(a) X_{r,s})(g_t^T(b) X_{r,t})(g_s^T(m+1) X_{r,s})(g_t^T(m+1) X_{r,t}). \end{aligned}$$

It can be shown that $G$ is a strongly consistent estimator of the covariance matrix of $j_k(Q)$ restricted to the tangent space of $j_k((\mathbb{R}P^m)^q)$ at $j_k(\mu)$, with respect to the orthobasis determined by the eigenvectors $g_s(a), s = 1, \ldots, q, a = 1, \ldots, m$.

Let $D_s = (g_s(1), \ldots, g_s(m))$, $s = 1, \ldots, q$, so that $D_s$ is an $(m+1) \times m$ matrix. If $\mu = ([\gamma_1], \ldots, [\gamma_q])$, where $\gamma_s$, $s = 1, \ldots, q$, are unit column vectors in $\mathbb{R}^{m+1}$, we define a Hotelling $T^2$ type-statistic,

$$(4.11) \qquad T^2(Y; Q) = n(\gamma_1^T D_1, \ldots, \gamma_q^T D_q) G^{-1} (\gamma_1^T D_1, \ldots, \gamma_q^T D_q)^T.$$

Note that $(\gamma_1^T D_1, \ldots, \gamma_q^T D_q)$ is a row vector. A proof of Theorem 4.1 regarding the asymptotic distribution of $T(Y; Q)$ is available in a technical report of Mardia and Patrangenaru [19]; the report also examines Fréchet's intrinsic mean.

THEOREM 4.1. *Assume* $(Y_r), r = 1, \ldots, n$, *is a random sample from a probability measure $Q$ on $(\mathbb{R}P^m)^q$ as above and for $s = 1, \ldots, q$, let $\lambda_s(a)$ and $\gamma_s(a)$ be the eigenvalues in increasing order and corresponding unit eigenvectors of $E[X_{1,s} X_{1,s}^T]$. If $\lambda_s(m+1) > 0$, $s = 1, \ldots, q$, are simple, then $T^2(Y; Q)$ in* (4.11) *converges weakly to $\chi_M^2$.*

The following result follows the technique of bootstrapping the distribution of extrinsic sample means on a manifold [3].

COROLLARY 4.1. *Let* $(Y_r), r = 1, \ldots, n$, *be a random sample from* $Q$ *on* $(\mathbb{R}P^m)^q$, *and let* $Y_r = ([X_{r,1}], \ldots, [X_{r,q}])$, $X_{r,j}^T X_{r,j} = 1, j = 1, \ldots, q$. *Assume that* $Q$ *has a nonzero absolutely continuous component. For a random resample* $(Y_1^*, \ldots, Y_n^*)$ *from* $(Y_1, \ldots, Y_n)$, *denote the eigenvalues of* $\frac{1}{n} \sum_{r=1}^{n} X_{r,s}^* X_{r,s}^{*T}$ *in increasing order by* $d_s^*(a), a = 1, \ldots, m+1$, *and the corresponding unit eigenvectors by* $g_s^*(a), a = 1, \ldots, m+1$. *Let* $G^*$ *be the matrix obtained from* $G$, *by substituting all the entries with* ∗*-entries. Then the bootstrap distribution function of the statistic*

$$
(4.12) \quad \begin{aligned} T^{*2}(Y^*, \hat{Q}_n) &= n(g_1^T(m+1)D_1^*, \ldots, g_q^T(m+1)D_q^*) \\ &\quad \times G^{*-1}(g_1^T(m+1)D_1^*, \ldots, g_q^T(m+1)D_q^*)^T \end{aligned}
$$

*approximates the true distribution of* $T^2(Y; Q)$, *given by* (4.11), *with an error of order* $O_p(n^{-2})$.

In imaging applications, we may be given a template (population) mean in terms of $k$ landmarks. The aim could then be to test for identification, that is, to assess if a given sample of images of a scene is from this population. Such a test can be based asymptotically on the test statistic given in (4.11). For small samples we need to use the bootstrapped statistic (4.12).

We compare extrinsic means of projective shapes with respect to the embedding $j_k$ of $P\Sigma_m^k = (\mathbb{R}P^m)^q$ in $(S(m+1))^q$. Using our axial representation, we reduce this problem to axial statistics. Since for $b = 1, 2$, $\mu_b \in (\mathbb{R}P^m)^q$, we set in axial representation $\mu_b = (\mu_{b,1}, \ldots, \mu_{b,q})$. Without loss of generality, we may assume that for each $j$ the angle between $\mu_{1,j}$ and $\mu_{2,j}$ is $\frac{\pi}{2}$ or less, and we consider the unique rotation $\rho_{m,j} \in SO(m+1)$ such that $\rho_{m,j}(\mu_{1,j}) = \mu_{2,j}$ and the restriction of $\rho_{m,j}$ to the orthocomplement of the plane determined by $\mu_{1,j}$ and $\mu_{2,j}$ in $\mathbb{R}^{m+1}$ is the identity.

The equality $\mu_1 = \mu_2$ is equivalent to $\rho_{m,s} = I_{m+1}, s = 1, \ldots, q$, where $I_{m+1}$ is an $(m+1) \times (m+1)$ identity matrix. Assume $(Y_{1,r}), r = 1, \ldots, n_1$, $(Y_{2,t}), t = 1, \ldots, n_2$, are random samples from $Q_1, Q_2$, respectively. A consistent estimator of the Lie group valued parameter $\rho_m = (\rho_{m,s}, s = 1, \ldots, q)$ is $r_m = (r_{m,s}, s = 1, \ldots, q)$, where, for each $s = 1, \ldots, q$, $r_{m,s} \in SO(m+1)$ is the unique rotation defined as above. This rotation brings the extrinsic sample means (mean directions) in coincidence, that is, superimposes $m_{1,s}$ onto $m_{2,s}$. Here $m_{a,s}$ is the unit eigenvector of $\sum_{s=1}^{q} X_{r,as} X_{r,as}^T$, where $Y_{a,r} = ([X_{r,a,1}], \ldots, [X_{r,a,q}])$, $r = 1, \ldots, m_a, a = 1, 2$.

*Case for $m=2$.* A particular case of practical interest is when $m=2$. We write here, for this particular case, $\rho$ for $\rho_2$ and $r$ for $r_2$. Here we will consider only the subcase $k=5$, for which we give an application in the next section. To test the equality $\mu_1 = \mu_2$ amounts to testing

$$\rho = id_{\mathbb{R}P^2},$$

where $id_{\mathbb{R}P^2}$ is the identity map of $\mathbb{R}P^2$ in the group of isometries $I(\mathbb{R}P^2)$ of $\mathbb{R}P^2$. Any isometry $\rho$ of $\mathbb{R}P^2$ can be represented in a unique way by a rotation $T \in SO(3)$ of the angle $\theta$, with $0 \leq \theta < \pi, \rho([x]) = [TX]$. Note that a rotation $T$ of the angle $\theta = \pi$ acts as the identity $id_{\mathbb{R}P^2}$. If $\rho \neq id_{\mathbb{R}P^2}$ and $T \in SO(3)$ represents $\rho$, there is an orthonormal basis $V_1, V_2, V_3$ of $\mathbb{R}^3$ such that $T(V_3) = V_3$. We set

$$H(\rho) = [(V_1 \cdot T(V_1), (V_1 \times T(V_1)))^T],$$

for $\rho \neq id_{\mathbb{R}P^2}$ and $H(id_{\mathbb{R}P^2}) = [1:0:0:0]$. The map $H: I(\mathbb{R}P^2) \to \mathbb{R}P^3$ is a well-defined isomorphism from $I(\mathbb{R}P^2)$ to the axial space $\mathbb{R}P^3$. Modulo the diffeomorphism $H$, the equality $\mu_1 = \mu_2$ amounts to $H(\rho) = [1:0:0:0]$. The distribution of the resulting consistent estimator $H(r)$ of $H(\rho)$ is essentially given in [1], Theorem 2.1. Assume neither $n_1$ nor $n_2$ is small compared with $n = n_1 + n_2$. Let $G(\rho)$ be the affine coordinates of $H(\rho)$; if $H(\rho) = [H_0(\rho):H_1(\rho):H_2(\rho):H_3(\rho)]$, then $G(\rho) = (G_1(\rho), G_2(\rho), G_3(\rho))$, with $G_a(\rho) = H_a(\rho)/H_0(\rho), a = 1, 2, 3$. Using equation (5.13) of [1], page 489, it turns out that $n^{1/2}(r - \rho)$ has asymptotically a trivariate Gaussian distribution and is independent of $n$. Then by the delta method of Cramér (see, e.g., [7], page 45), it follows that, under the null hypothesis, if there are two constants $c_1, c_2$, such that $n_b/n \to c_b$, for $b = 1, 2$, then $n^{1/2}\{G(r) - G(\rho)\}$ has asymptotically a trivariate Gaussian distribution which is independent of $n$. Consequently, if we consider the resamples under the empirical distribution $n^{1/2}\{G(r^*) - G(r)\}$ by a nonpivotal bootstrap, then this distribution will have asymptotically the same distribution as that of $n^{1/2}\{G(r) - G(\rho)\}$.

*Concentrated data case.* For each projective coordinate from concentrated projective shape data, we may select only one representative on the sphere. Therefore, we may use the directional representation on a product of $q$ copies of $S^m$.

In this case, if $Q$ is a concentrated distribution on $(S^m)^q$, let $\mu_D = (\frac{\mu_1}{\|\mu_1\|}, \ldots, \frac{\mu_q}{\|\mu_q\|})$ be the mean multivariate direction and $\bar{y}_D = (\frac{\bar{y}^1}{\|\bar{y}^1\|}, \ldots, \frac{\bar{y}^q}{\|\bar{y}^q\|})$ be the corresponding sample mean corresponding to a random sample $(y_1, \ldots, y_n)$, $y_j = (y_j^1, \ldots, y_j^q) \in (S^m)^q$. The asymptotic distribution of $\bar{y}_D$ can be described in terms of an orthonormal frame field $(e_{1,1}(y^1), \ldots, e_{1,m}(y^1), \ldots, e_{q,1}(y^q), \ldots, e_{q,m}(y^q))$ defined around $\mu_D$; here $y = (y^1, \ldots, y^q) \in (S^m)^q$ and for each $a = 1, \ldots, q, e_{a,i}(y^a)^T e_{a,j}(y^a) = \delta_{ij}, i, j = 1,$

where $e_{a,i}(y^a) \in \mathbb{R}^{m+1}$. Let $G(y)$ be the $M$ by $M$ matrix made of $q$ by $q$ square matrices of size $m, G_{ab}(y)$, where

$$G_{ab}(y) = n^{-1}(\|\bar{y}^a\|\|\bar{y}^b\|)^{-1}\left(\sum_{r=1}^{n}\left(e_{a,i}\left(\frac{\bar{y}^a}{\|\bar{y}^a\|}\right)^T y_r^a\right)\left(e_{b,j}\left(\frac{\bar{y}^b}{\|\bar{y}^b\|}\right)^T y_r^b\right)\right),$$

$$i, j = 1, \ldots, m.$$

We studentize the tangential component of the difference between the sample and population directional means and obtain the following result.

THEOREM 4.2. *Let $\bar{y}_D - \mu_D = \sum_{a=1}^{q}\sum_{j=1}^{m} d_{a,j} e_{a,j}(\frac{\mu^j}{\|\mu^j\|}) + \nu$ be the decomposition of the difference between the directional sample mean and directional mean into its tangential and normal components. If $d_a = (d_{a,j})^T$, $j = 1, \ldots, m$, and $d = (d_1, \ldots, d_q)^T$, then $T^2(Y, Q, \mu_D) = n \cdot d^T G(y)^{-1} d$ converges weakly to $\chi^2_M$.*

REMARK 4.2. For concentrated data, Proposition 3.1 and Theorem 4.2 show that asymptotically the squared norm of the Studentized sample mean and of the Studentized extrinsic sample mean vector both have a $\chi^2_M$ distribution, where $M$ is the dimension of the projective shape space.

COROLLARY 4.2. *Assume $Y_1^*, \ldots, Y_n^*$ is a random resample with replacement of a random sample $(Y_1, \ldots, Y_n)$ from a probability measure $Q$ on $(S^m)^q$ such that $\mu_D$ exists and $Q$ has a unique absolutely continuous component with respect to the volume measure of $(S^m)^q$. Let $\hat{Q}_n$ be the empirical distribution, and let $d^*, G(y^*)$ be the corresponding bootstrap analogs of $d$ and $G(y)$ obtained by substituting $Y_1^*, \ldots, Y_n^*$ for $Y_1, \ldots, Y_n$, $\bar{Y}$ for $\mu$, and $\bar{Y}^*$ for $\bar{Y}$. Then the distribution function of $T^2(Y, Q, \mu_D)$ can be approximated by the bootstrap distribution of $T^2(Y^*, \hat{Q}_n, \bar{Y}_D) = nd^{*T}G(Y^*)^{-1}d^*$ with a coverage error of order $O_p(n^{-2})$.*

COROLLARY 4.3. *A $(1-\alpha)100\%$ bootstrap confidence region for $\mu_D$ is given by the following:*

(a) $R_\alpha(Y) = \{\mu \in (S^m)^q | T^2(Y, Q, \mu) \leq T_\alpha^{*2}\}$, *where $T_\alpha^{*2}$ is the $(1-\alpha)100\%$ percentile of the bootstrap distribution $T^2(Y^*, \hat{Q}_n, \bar{Y}_D)$.*

(b) $S_\alpha(Y) = \{\mu = (\mu_1, \ldots, \mu_q) \in (S^m)^q | T_j^2(Y, Q, \mu_j) \leq T_{j,\alpha/q}^{*2}, j = 1, \ldots, q\}$, *where $T_{j,\alpha/q}^{*2}$ is the $(1-\frac{\alpha}{q})100\%$ percentile of the bootstrap distribution $T_j^2(Y^*, \hat{Q}_n, \bar{Y}_D^j)$ corresponding to the $j$th directional component only.*

While Corollary 4.3(a) is useful when $M = 1$ or $M = 2$ (see Example 5.1), for larger values of $M$ the computations are very intensive. To decrease the computational complexity, one may use Corollary 4.3(b), which is based on a Bonferroni type of argument and gives the confidence region as a Cartesian product of confidence regions for the directional components of $\mu_D$. This is used in Example 5.3 below.

**5. Applications.** To illustrate our methodology, we consider machine vision applications involving data extracted from photographs. In the first two examples we use architectural features of two buildings, which form the preliminary stage of object recognition. In the third example we consider a problem in face recognition where two different views are available. We have assumed here that the coordinates of points are already recorded from a digital image; image processing software often has built-in procedures to extract coordinates of landmarks. There are various algorithms in practice to carry out this task (see, e.g., [12]). Our main aim is to show how we can use the one-sample and the two-sample tests when the underlying hypotheses are plausible. We show the type of computations required and indicate how the approximations work. We also examine how the parametric and the nonparametric tests perform. In Example 5.1 the problem is illustrated for $m=1, k=4$ so $q=1$ and involves circular statistics. In Example 5.2 $m=2$, $k=5$, whereas in Example 5.3 $m=2, k=6$, so $q=1$ and $q=2$, respectively, so that Example 5.2 involves the univariate spherical statistics, whereas Example 5.3 involves bivariate spherical statistics. These ideas can be extended to problems in machine vision such as in identification, classification and so on. Note that in these illustrative examples the number of landmarks $k$ and the sample size $n$ happened to be small, but the methods are applicable to any values of $k$ and $n$.

EXAMPLE 5.1. In this example we selected randomly five photographic views of a building (the Education Building) from a large database of buildings at the University of Leeds. One view is shown in Figure 5. It can be seen that the windows are equi-spaced. We selected four landmarks which are four central points of four consecutive windows (see Figure 5). The observed values of these landmarks are given in Table 1. Note that here $m=1, k=4, q=1$ and $n=5$.

We fix the projective frame $\pi = ([x_1:1], [x_2:1], [x_3:1])$ and determine the cross-ratio $c$ and projective coordinate of $p^\pi$ of $p = [x_4:1]$ given by (2.15). After doubling the angles, for each view, we get a direction $\theta$. These values are shown in Table 1, together with the coordinates.

If the landmarks are equidistant, their cross-ratio is $c = 4/3$ and the corresponding direction is $\theta_0 = 1.287$ rad. Therefore, testing the hypothesis for projective equidistance is equivalent to the problem

$$H_0 : \theta = \theta_0 \quad \text{vs.} \quad H_1 : \theta \neq \theta_0.$$

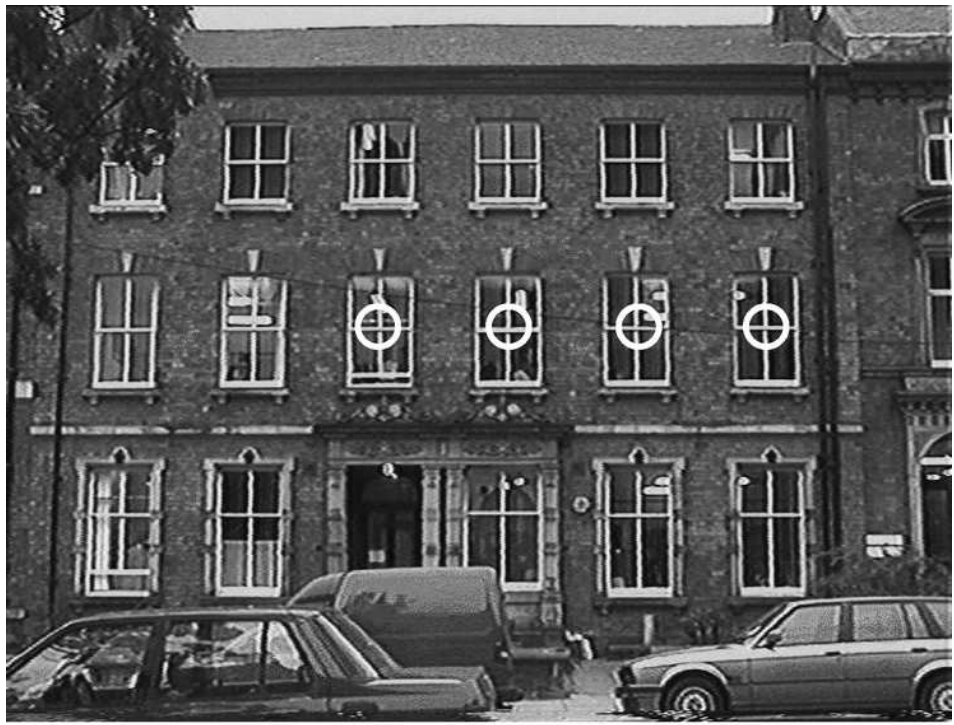

Education building: view 1

Fig. 5. *One view of the Education Building (University of Leeds) with four collinear points on four consecutive windows marked by white rings at the centre of the windows of the first floor.*

Table 1
*Horizontal coordinates of collinear points on the Education Building at 5 different views*

| View | $x_1$ | $x_2$ | $x_3$ | $x_4$ | $c$ | $\phi$ | $\theta$ |
|---|---|---|---|---|---|---|---|
| 1 | 22.90 | 35.7 | 48.3 | 61.10 | 1.340 | 0.641 | 1.282 |
| 2 | 23.10 | 29.1 | 35.5 | 42.50 | 1.338 | 0.642 | 1.284 |
| 3 | 41.40 | 44.3 | 47.3 | 50.70 | 1.353 | 0.636 | 1.273 |
| 4 | 39.00 | 47.0 | 53.9 | 60.00 | 1.337 | 0.642 | 1.285 |
| 5 | 42.25 | 46.9 | 50.5 | 53.85 | 1.373 | 0.629 | 1.259 |

*A parametric approach.* Under the assumption of the von Mises distribution, as in [11], we consider the Watson–Williams [28] test statistic

$$F^{(1)} = (n-1)(R - (\cos\theta_0, \sin\theta_0)\mathbf{R})/(n-R),$$

where $R$ is the length of the resultant column vector $\mathbf{R} = \sum_{r=1}^{n}(\cos\theta_r, \sin\theta_r)^T$. Under $H_0$, for concentrated data $F^{(1)}$ is approximately distributed as $F_{1,n-1}$. For this data $F^{(1)}$ was found to be 2.826, which has $p$-value 0.168. Hence, we fail to reject the null hypothesis at the 5% significance level.

Indeed, here $\bar{R} = R/n = 0.99994$, so that data is highly concentrated. Hence, we can assume the tangent approximation of Section 3 is valid. Under this assumption the one sample $t$-test yields the $p$-value 0.170. So again we fail to reject $H_0$ at the 5% significance level.

*A nonparametric approach.* Assume now that the population distribution $Q$ is arbitrary and has a mean direction $\mu_D = \exp(i\theta_D)$. We consider the hypotheses

$$H_0: \theta_D = \theta_0 \quad \text{vs.} \quad H_1: \theta_D \neq \theta_0.$$

Since the sample size $n = 5$ is very small, we base our $p$-value on Corollary 4.1. Using 5000 resample values of $T^2(Y^*, \hat{Q}_n, \bar{Y}_D)$, we found the $p$-value to be 0.201. Thus, we again fail to reject $H_0$ at the 5% significance level. In conclusion, we fail to reject the equidistance hypothesis at the 5% level using either test.

EXAMPLE 5.2. In this example we illustrate the two sample tests for $m = 2$ and $k = 5$ so that $q = 1$. Again, we have used the Leeds University Buildings database. In addition to the Education Building used in the previous example, we now consider an additional building—the Careers Building. Two groups of identically positioned noncollinear landmarks $A_1, A_2, A_3, A_4, A_5$ were marked on five frontal photographs of the Education Building and four of the Careers Building, so that $n_1 = 4$ and $n_2 = 5$. One of the buildings with the landmarks is shown in Figure 6.

We obtain the spherical coordinates of the landmarks for the two samples following the calculations similar to those in Example 2.1, and these are given in Table 2. Assessing this part of architectural similarity between the two buildings is equivalent to performing a two-sample test for means. It is clear from our images that the architectural style of the windows based on these landmarks is very similar. We will also show that the Hotelling $T^2$ test based on projective invariants leads to a contradictory result, which indicates that we should prefer the use of the spherical projective coordinates.

*A parametric test.* For the Education Building and for the Careers Building, we find using Table 2 that the mean resultant lengths are 0.9997 and 0.9979, respectively. Hence, the data are highly concentrated in projective

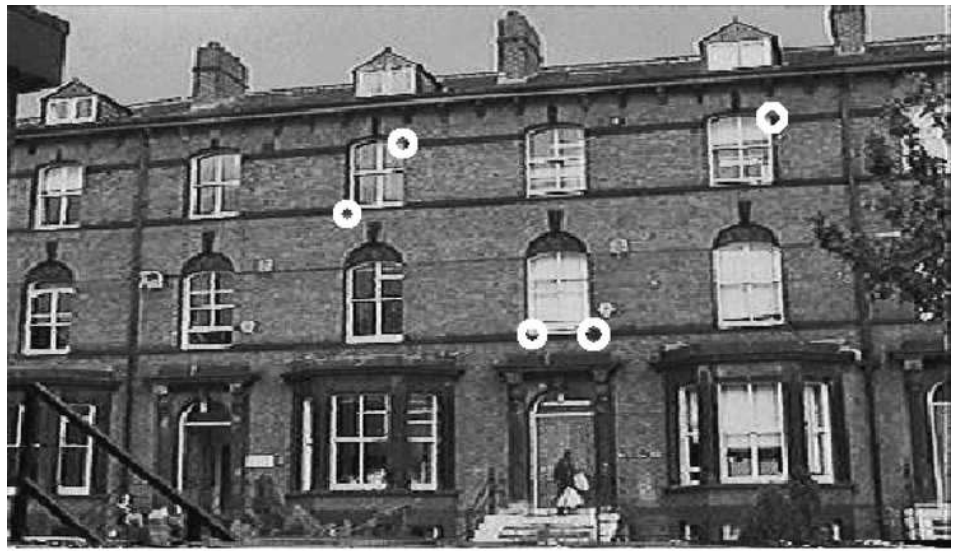

Careers building: view 1

Fig. 6. *One view of the Careers Buildings (University of Leeds) with five landmarks.*

shape space. Thus, we use Hotelling's $T^2$ test in the tangent space given by (3.5). The combined mean direction and the mean resultant lengths are

$$(0.7980, 0.5722, 0.1892), 0.9988,$$

respectively. We find the value of $F$ for Hotelling's $T^2$ test is 2.6075, which is $F_{2,6}$. In fact, $Pr(F_{2,6} > 2.6075) \simeq 0.225$ so that we fail to reject the null hypothesis.

*A nonparametric test.* We selected the projective frame $\pi = ([A_1:1], [A_2:1], [A_3:1], [A_4:1])$ and determined the coordinates of the views in the sample, using a spherical representation; these spherical coordinates are displayed in Table 2.

Table 2
*Spherical coordinates for the Education Building (five views) and the Careers Building (four views)*

| View | Education Building | | | View | Careers Building | | |
|---|---|---|---|---|---|---|---|
| 1 | 0.8142 | 0.5547 | 0.1718 | 1 | 0.7859 | 0.5768 | 0.2228 |
| 2 | 0.8038 | 0.5610 | 0.1977 | 2 | 0.8170 | 0.5712 | 0.0791 |
| 3 | 0.8067 | 0.5591 | 0.1917 | 3 | 0.7639 | 0.6041 | 0.2268 |
| 4 | 0.8150 | 0.5513 | 0.1787 | 4 | 0.7893 | 0.5766 | 0.2110 |
| 5 | 0.7773 | 0.5890 | 0.2211 | | | | |

Here the extrinsic sample mean projective shapes of views from the Education Building and Careers Building are given in the spherical representation by $\bar{Y}_{1,E} = [0.8037 : 0.5632 : 0.1922]$ and $\bar{Y}_{2,E} = [0.7907 : 0.5834 : 0.1855]$, respectively. Now consider the problem of estimating the distribution of the axis $H(r)$ defined in Section 4. Since the smaller sample size is 4, and the eigenanalysis has to be repeated for each resample, we limited ourselves to 250 pseudorandom resamples and determined the corresponding nonpivotal bootstrap distribution of $G(r^*)$. The corresponding distribution of $3G(r^*)$ is displayed in Figure 7, which indicates that the sample mean of $G(r)$ is close to $(0, 0, 0)$.

The rotation that brings $\bar{Y}_{1,E}$ in coincidence with $\bar{Y}_{2,E}$ is identified with a 4-D axis (see Section 4), which turns out to be

$$H(r) = [0.9997 : -0.0077 : 0.0029 : 0.0231],$$

where we have used the dot product and cross-product of $\bar{Y}_{1,E}$ and $\bar{Y}_{2,E}$. We determined the coordinates of the distribution of $3(G(r^*) - G(r))$ and, for this distribution, we successively sorted and trimmed the distribution of $3G(r^*) = \{3G_1(r^*), 3G_2(r^*), 3G_3(r^*)\}$ and obtained the following 93% simultaneous bootstrap confidence intervals: $[-4.36, 3.02]$ for $3G_1(r^*)$, $[-3.59, 2.67]$ for $3G_2(r^*)$, $[-2.70, 3.40]$ for $3G_3(r^*)$. This analysis shows that $(0, 0, 0)$ is in the 93% percentile confidence region, that is, the identity is in the corresponding 93% bootstrap confidence region for $\rho_2$ on $SO(3)$. Therefore, we fail to reject $\mu_1 = \mu_2$ at significance level $\alpha = 0.07$.

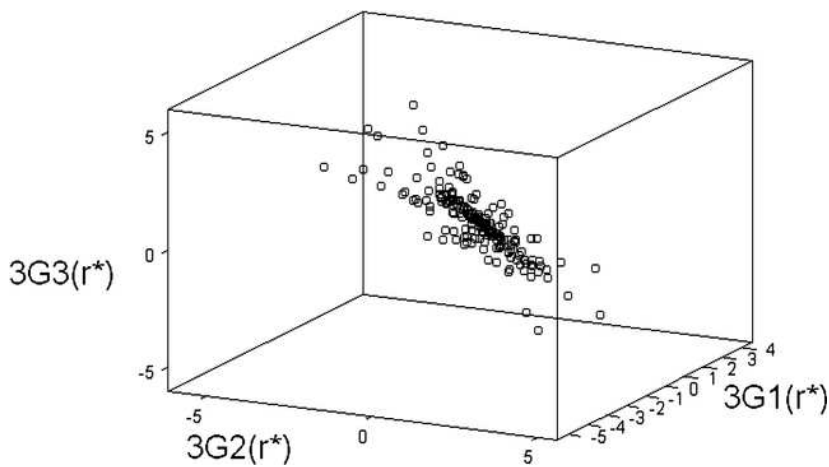

FIG. 7. *Affine view of the bootstrapped distribution of size* 250 *of the nonpivotal vector* $(3G1(r^*), 3G2(r^*), 3G3(r^*))$ *used in the comparison of the projective shapes of five landmarks on the Education and Careers Buildings.*

*A test based on invariants.* Using Table 2, the projective invariants $\iota_1$ and $\iota_2$, defined at (2.16), are given below for the two buildings, respectively:

Education Building: (4.739, 3.229), (4.068, 2.838), (4.208, 2.917), (4.561, 3.083), (3.516, 2.664),

Careers Building: (3.527, 2.588), (10.325, 7.219), (3.369, 2.664), (3.741, 2.733).

We find that the observed value of $F$ for the Hotelling $T^2$ test is 12.22 and $P(F_{2,4} > 12.22) = 0.0077$. Thus, we reject the hypothesis of similarity between the two buildings. This conclusion is quite different than the one we get using projective spherical shape coordinates. Indeed, this aspect of the architecture is so similar that we should be accepting the hypothesis. The difference is explained by the fact that the data is not normal, and the test based on invariants is sensitive to departures from normality.

EXAMPLE 5.3. We now apply the method to a face recognition problem. Figure 8(a) shows the seven frontal views ($n_1 = 7$) of the same person (an actor posing in different disguises) and Figure 8(b) shows his seven side views ($n_2 = 7$). We recorded six landmarks (four corners of the eyes, "canthus," and two end points of the lips, "mouth edge points"). Using the four eye-corner landmarks as the projective frame, the Cartesian landmarks were converted into the directional representation (bivariate spherical), leading to the spherical projective coordinates $x_1$ and $x_2$ in the same way as in Example 2.1. This data is displayed in Table 3.

*A parametric test.* First, it can be seen that, for frontal views, the mean resultant lengths of $x_1$ and $x_2$ are 0.9995, 0.9955, respectively, whereas for the side views, the mean resultant lengths are 0.9995, 0.9996, respectively. These values imply that the data is highly concentrated. For the combined data, the respective mean resultant lengths are 0.9995, 0.9996. Thus, we could use the tangent space to test the hypothesis that the two means are equal (Section 3). We find that the combined mean directions are given by

$$\hat{\mu}_1^T = (0.6889, 0.6735, 0.2681), \qquad \hat{\mu}_2^T = (0.7015, 0.6874, 0.1882).$$

We calculated the tangent coordinates using (3.4), which, from (3.5), leads to the value of $F = 0.8269$; this has an $F$-distribution with degrees of freedom 4 and 5. Since $Pr(F_{4,9} > 0.8269) = 0.5402$, there is strong evidence that these frontal and side views are of the same person.

*A nonparametric test.* The 95% pivotal bootstrap confidence region in the mean bivariate spherical direction based on Corollary 4.3(b), using 1500 random resamples was found to be $T^{*2}_{1,0.025} = 3.21$ and $T^{*2}_{2,0.025} = 1.84$ and the statistics $T_1^2(\bar{y}_D, y), T_2^2(\bar{y}_D, y)$ obtained, using the seven side views, were 1.54 and 1.33, respectively. It is seen that the bivariate mean projective shape

(a)

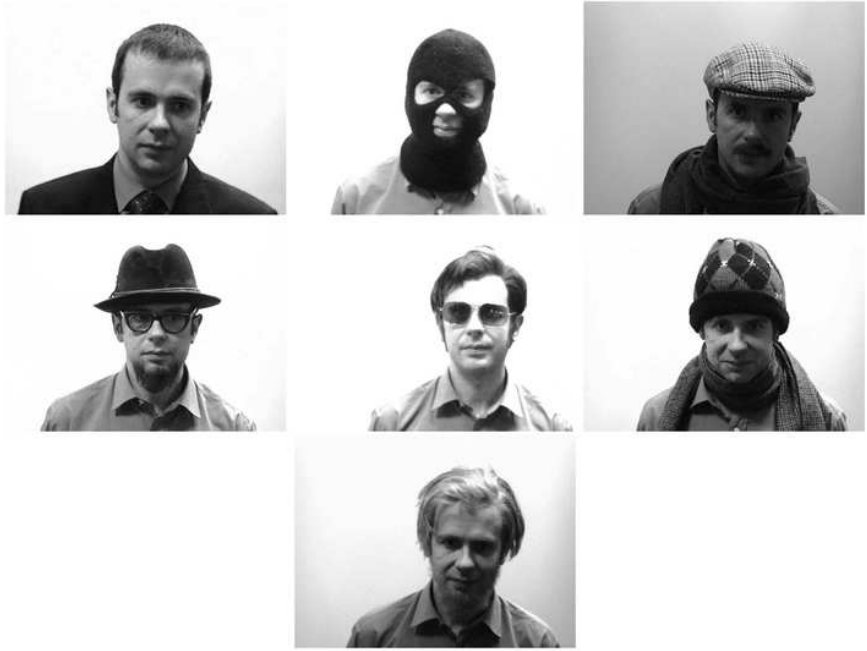

(b)

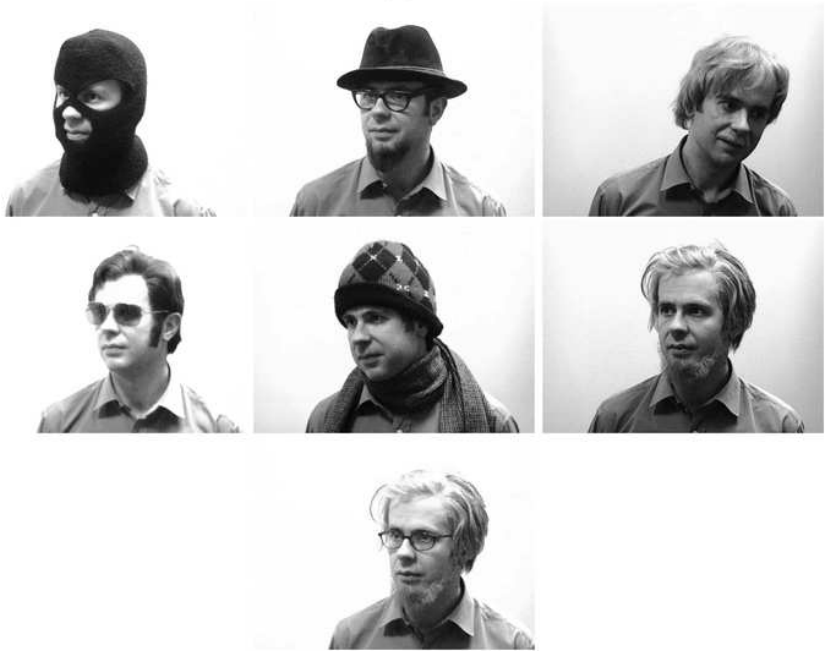

Fig. 8. *Photographs of an actor posing in different disguises.* (a) *Frontal views.* (b) *Side views.*

TABLE 3
*Bivariate spherical coordinates $x_1 = (x_1^1, y_1^1, z_1^1)$ and $x_2 = (x_2^2, y_2^2, z_2^2)$ corresponding to the 14 views (1–7 frontal, 8–14 side) of the actor*

| View | $x_1^1$ | $y_1^1$ | $z_1^1$ | $x_2^2$ | $y_2^2$ | $z_2^2$ |
|---|---|---|---|---|---|---|
| 1 | 0.700780 | 0.657783 | 0.276096 | 0.708981 | 0.676761 | 0.198344 |
| 2 | 0.685337 | 0.675546 | 0.271939 | 0.697420 | 0.691293 | 0.188996 |
| 3 | 0.688405 | 0.650635 | 0.320581 | 0.709839 | 0.669692 | 0.218266 |
| 4 | 0.690658 | 0.673332 | 0.263846 | 0.706231 | 0.681266 | 0.192651 |
| 5 | 0.691832 | 0.668204 | 0.273626 | 0.700515 | 0.685421 | 0.198688 |
| 6 | 0.688246 | 0.667378 | 0.284470 | 0.703869 | 0.680057 | 0.205158 |
| 7 | 0.681884 | 0.685515 | 0.255155 | 0.692303 | 0.697768 | 0.183948 |
| 8 | 0.679369 | 0.669555 | 0.300255 | 0.694591 | 0.683580 | 0.224191 |
| 9 | 0.686636 | 0.687718 | 0.235742 | 0.698689 | 0.696648 | 0.162835 |
| 10 | 0.684002 | 0.685272 | 0.250087 | 0.694651 | 0.701190 | 0.160603 |
| 11 | 0.667353 | 0.699274 | 0.256235 | 0.679292 | 0.713641 | 0.171112 |
| 12 | 0.717523 | 0.665701 | 0.204944 | 0.726679 | 0.673698 | 0.134419 |
| 13 | 0.694639 | 0.669937 | 0.262035 | 0.705996 | 0.686475 | 0.174132 |
| 14 | 0.723910 | 0.649280 | 0.233215 | 0.735734 | 0.656621 | 0.165965 |

corresponding to the seven side views falls in the 95% pivotal bootstrap confidence region for the mean projective shape from frontal views. Thus, both methodologies lead to the same conclusion.

REMARK 5.1. We have used an edge registration method in developing the projective coordinates which are points in a curved space. The underlying projective shape space has the features of Kendall's shape space, since it is a symmetric space. For large samples, the choice of edge registration will have no influence on the analysis. For a similar discussion on the Bookstein shape coordinates, see [5], page 30.

**6. A strategy for general shape analysis.** We now give a unified strategy underlying the three statistical shape spaces: similarity, affine and projective, leading to general statistical shape analysis. Each of these spaces is a space of orbits (orbifolds) of group actions on a finite set of points on a manifold. In general shape analysis, the group actions of interest are the following:

1. In the case of similarity shape the group of direct similarities or, more generally, the group $CO(m) \ltimes \mathbb{R}^m$ of conformal linear maps of $\mathbb{R}^m$, and the manifold is $\mathbb{R}^m$.
2. In the case of affine shape, the group $GL(m, \mathbb{R}) \ltimes \mathbb{R}^m$ of affine transformations of $\mathbb{R}^m$ and the manifold is $\mathbb{R}^m$.
3. In the case of projective shape, the group $PGL(m, \mathbb{R})$ of projective transformations of $\mathbb{R}P^m$ and the manifold is $\mathbb{R}P^m$, or the pseudo-group action of $PGL(m, \mathbb{R})$ on $\mathbb{R}^m$, regarded as the subset $\mathbb{R}P^m \setminus \mathbb{R}P^{m-1}$ of $\mathbb{R}P^m$.

Note that the groups are hierarchically ordered in these three types of shape.

Let $k$ be a fixed positive integer, and let $\mathcal{C}_m^k$ be the set of all configurations of $k$ labelled points in $\mathbb{R}^m$ (or in $\mathbb{R}P^m$). Assume $G$ is a group of transformations of $\mathbb{R}^m$ (or of $\mathbb{R}P^m$) that acts on the left on $\mathcal{C}_m^k$, via

$$g((x_1,\ldots,x_k)) = (g(x_1),\ldots,g(x_k)), \qquad g \in G, x = (x_1,\ldots,x_k) \in \mathbb{R}^m.$$

The *orbit* $G(x)$ of such a configuration $x = (x_1,\ldots,x_k) \in \mathcal{C}_m^k$ is defined by $G(x) =: \{g(x)|g \in G\}$. The *full $G$-shape space*, $G\Sigma_m^k$, is the set of all $G$-orbits, $G\Sigma_m^k = \{G(x), x = (x_1,\ldots,x_k) \in \mathcal{C}_m^k\}$. Note that the similarity planar shape space $\Sigma_2^k$ (see [14]) is not a full shape space, since the orbit of $x = (0,\ldots,0)$ is removed from the full $CO(2) \ltimes \mathbb{R}^2$-shape space. If this singular orbit is removed, the corresponding shape space has the structure of a manifold, namely, the complex projective space $\mathbb{C}P^{k-2} \cong \Sigma_m^k$.

The strategy used for similarity shape or affine shape can be extended to the general context by taking the following sequence of steps for any group action on $\mathcal{C}_m^k$ whose orbits are closed:

1. Identify $\mathcal{C}_m^k$ with $\mathbb{R}^{km}$ via $(x_1,\ldots,x_k) \to (x_1^1,\ldots,x_1^m,\ldots,x_k^1,\ldots,x_k^m)$, $x_j = (x_j^1,\ldots,x_j^m), j = 1,\ldots,k$, and consider the topology on $\mathcal{C}_m^k$ inherited from the Euclidean topology of $\mathbb{R}^{km}$.
2. Consider the quotient topology on the orbit space $G\Sigma_m^k$ and let $\pi: \mathcal{C}_m^k \to G\Sigma_m^k$ be the quotient map. A subset $\mathcal{U}$ of $G\Sigma_m^k$ is open if $\pi^{-1}(\mathcal{U})$ is open in $\mathcal{C}_m^k$. Recall that a subset $\mathcal{V}$ is *generic* if it is open and dense in the quotient topology. Note that if $\pi: \mathcal{C}_m^k \to G\Sigma_m^k$ is the quotient map, then $\mathcal{V}$ is generic if $\pi^{-1}(\mathcal{V})$ is a generic subset in $\mathcal{C}_m^k$. Consider a generic subset $\mathcal{V} \subseteq G\Sigma_m^k$ that has a homogeneous structure (see [27]) or even a structure of symmetric space.
3. Whenever possible, find generic shape spaces that admit homogeneous structures. In such a situation find an equivariant embedding of $\mathcal{V}$ into a Euclidean space which yields easily computable extrinsic sample means.
4. Derive the distributions for marginal distributions on $\mathcal{V}$ resulting from noise at landmark locations, and if these distributions are intractable, approximate them with simpler distributions.
5. Determine asymptotic distributions of Fréchet sample means on $\mathcal{V}$ and, in particular, of extrinsic sample means. Then use associated statistics to design large sample confidence regions for population extrinsic means. For small samples derive corresponding bootstrap distributions for sample means and confidence regions.

In this paper we introduced the projective shape space of configurations of points in general position following this strategy. In this case, given that the projective transformations in $\mathbb{R}^m$ form only a pseudogroup, it was more convenient to regard $\mathbb{R}^m$ as an open affine set of $\mathbb{R}P^m$, and to consider points in $\mathbb{R}P^m$ and the group $G = PGL(m)$. We considered the generic set $\mathcal{V} \subseteq G\Sigma(\mathbb{R}P^m)^k$ of projective shapes of configurations of $k$ points in $\mathbb{R}P^m$ with $k \geq m+2$ in general position, for which the first $m+2$ of these points form a projective frame in $\mathbb{R}P^m$. This generic set $\mathcal{V}$ was called $P\Sigma_m^k$ and we showed that $P\Sigma_m^k$ is a manifold diffeomorphic with $(\mathbb{R}P^m)^q$, thus having a structure of homogeneous space. We embedded this manifold in a space of matrices, and computed extrinsic sample means and their asymptotic distributions and derived bootstrap results to deal with small sample sizes. In the case of linear projective shapes, we also approximated distributions on the space of projective shapes, resulting from noise at landmark locations, by distributions that are easier to handle, thus completing the general program presented above. In the case of planar similarity shape, for $k > 2$, the group of similarities acts freely on the space of configurations $\mathcal{C}_2^k$. In general, if the restriction of the action of the group $G$ on $\mathcal{C}_m^k$ to a generic subset of orbits is free, the dimension of each orbit is equal to the manifold dimension of the Lie group $G$. In this case one may locally select a submanifold of $\mathcal{C}_m^k$ that is transverse to all the orbits, that is, $\dim(G\Sigma_m^k) = km - \dim G$. Table 4 gives the dimensions of similarity, affine and projective shape spaces. Note that the number of degrees of freedom of the chi-square distributions needed for confidence regions of Fréchet means is equal to the dimension of $G\Sigma_m^k$.

REMARK 6.1. Although less studied in statistics, projective shape analysis is the most relevant in image analysis, since the pinhole camera principle is based on central projections. Affine shape analysis and similarity shape analysis are valid in image analysis only when such central projections can be approximated with parallel projections, or even orthogonal projections.

REMARK 6.2. Comparison of projective shapes is made easier due to homogeneity of the projective shape space. Recall that a space $\mathcal{M}$ is homogeneous if there is a Lie group $G$ of transformations such that, for any points $x_1, x_2 \in \mathcal{M}$, there is a $g \in G$ with $g(x_1) = x_2$. In this way we may define a map from $\mathcal{M} \times \mathcal{M}$ to $G$ [which is what we did in Example 5.2, where

TABLE 4
*The appropriate dimensions for different shape spaces with k points in a configuration*

| Shape type | Similarity | Affine | Projective |
|---|---|---|---|
| Group | $CO(m) \ltimes \mathbb{R}^m$ | $\text{Aff}(m, \mathbb{R})$ | $PGL(m, \mathbb{R})$ |
| Dimension | $\frac{m(m+1)}{2} + 1$ | $m(m+1)$ | $m(m+2)$ |
| Dimension of shape space | $mk - \frac{m(m+1)}{2} - 1$ | $m(k-m-1)$ | $m(k-m-2)$ |

34 K. V. MARDIA AND V. PATRANGENARU

$\mathcal{M} = P\Sigma_2^5$ and $G = SO(3)$] and the comparison of two means is transferred on $G$. This method of comparison of projective shapes can be used to compare means of two populations on an arbitrary Riemannian homogeneous manifold and, in particular, on a Grassmanian manifold. Indeed, recently a need for population means on Grassmannian manifolds has arisen from signal processing; see [26].

**Acknowledgments.** We are most grateful to an Associate Editor, the two anonymous referees, and to Fred Bookstein and John Kent for their comments on various revisions which have led to substantial improvements. Our thanks also go to our collaborators David Hogg for the Leeds Building images and to Zeda Ltd. for the actor's photographs. Thanks are also due to Paul McDonnell, Ian Moreton, Vysaul Nyirongo, Robert Page, Ray Pruett, Charles Taylor and Alistair Walder for their help.

## REFERENCES

[1] BERAN, R. and FISHER, N. I. (1998). Nonparametric comparison of mean directions or mean axes. *Ann. Statist.* **26** 472–493. MR1626051
[2] BHATTACHARYA, R. N. and PATRANGENARU, V. (2003). Large sample theory of intrinsic and extrinsic sample means on manifolds. I. *Ann. Statist.* **31** 1–29. MR1962498
[3] BHATTACHARYA, R. N. and PATRANGENARU, V. (2005). Large sample theory of intrinsic and extrinsic sample means on manifolds. II. *Ann. Statist.* **33** 1225–1259. MR1962498
[4] BOOKSTEIN, F. L. (1991). *Morphometric Tools for Landmark Data*. Cambridge Univ. Press. MR1469220
[5] DRYDEN, I. L. and MARDIA, K. V. (1998). *Statistical Shape Analysis*. Wiley, Chichester. MR1646114
[6] FAUGERAS, O. and LUONG, Q.-T. (2001). *The Geometry of Multiple Images*. MIT Press, Cambridge, MA. MR1823958
[7] FERGUSON, T. (1996). *A Course in Large Sample Theory*. Chapman and Hall, London. MR1699953
[8] FISHER, N. I., HALL, P., JING, B.-Y. and WOOD, A. T. A. (1996). Improved pivotal methods for constructing confidence regions with directional data. *J. Amer. Statist. Assoc.* **91** 1062–1070. MR1424607
[9] FRÉCHET, M. (1948). Les éléments aléatoires de nature quelconque dans un espace distancié. *Ann. Inst. H. Poincaré* **10** 215–310. MR27464
[10] GOODALL, C. R. and MARDIA, K. V. (1993). Multivariate aspects of shape theory. *Ann. Statist.* **21** 848–866. MR1232522
[11] GOODALL, C. R. and MARDIA, K. V. (1999). Projective shape analysis. *J. Comput. Graph. Statist.* **8** 143–168. MR1706381
[12] HARTLEY, R. and ZISSERMAN, A. (2000). *Multiple View Geometry in Computer Vision*. Cambridge Univ. Press. MR2059248
[13] HEYDEN, A. (1995). Geometry and algebra of multiple projective transformations. Ph.D. dissertation, Univ. Lund, Sweden. MR1420745
[14] KENDALL, D. G. (1984). Shape manifolds, Procrustean metrics and complex projective spaces. *Bull. London Math. Soc.* **16** 81–121. MR737237

[15] KENDALL, D. G., BARDEN, D., CARNE, T. K. and LE, H. (1999). *Shape and Shape Theory.* Wiley, New York. MR1891212

[16] KENT, J. T. (1992). New directions in shape analysis. In *The Art of Statistical Science. A Tribute to G. S. Watson* (K. V. Mardia, ed.) 115–127. Wiley, New York. MR1175661

[17] MARDIA, K. V. (1975). Statistics of directional data (with discussion). *J. Roy. Statist. Soc. Ser. B* **37** 349–393. MR402998

[18] MARDIA, K. V. and JUPP, P. E. (2000). *Directional Statistics.* Wiley, New York. MR1828667

[19] MARDIA, K. V. and PATRANGENARU, V. (2002). Directions and projective shapes. Technical Report No. 02/04, Dept. Statistics, Univ. Leeds.

[20] MAYBANK, S. J. and BEARDSLEY, P. A. (1994). Classification based on the cross ratio. *Applications of Invariance in Computer Vision. Lecture Notes in Comput. Sci.* **825** (J. L. Mundy, A. Zisserman and D. Forsyth, eds.) 433–472. Springer, Berlin.

[21] PRENTICE, M. J. (1984). A distribution-free method of interval estimation for unsigned directional data. *Biometrika* **71** 147–154. MR738335

[22] PRENTICE, M. J. and MARDIA, K. V. (1995). Shape changes in the plane for landmark data. *Ann. Statist.* **23** 1960–1974. MR1389860

[23] SINGH, H., HNIZDO, V. and DEMCHUK, E. (2002). Probabilistic model for two dependent circular variables. *Biometrika* **89** 719–723. MR1929175

[24] SPARR, G. (1992). Depth computations from polyhedral images. *Image and Vision Computing* **10** 683–688.

[25] SPIVAK, M. (1970). *A Comprehensive Introduction to Differential Geometry.* Publish or Perish, Boston.

[26] SRIVASTAVA, A. and KLASSEN, E. (2002). Monte Carlo extrinsic estimators of manifold-valued parameters. *IEEE Trans. Signal Process.* **50** 299–308.

[27] TRICERRI, F. and VANHECKE, L. (1983). *Homogeneous Structures on Riemannian Manifolds.* Cambridge Univ. Press. MR712664

[28] WATSON, G. S. and WILLIAMS, E. (1956). On the construction of significance tests on the circle and the sphere. *Biometrika* **43** 344–352. MR82243

[29] ZIEZOLD, H. (1977). On expected figures and a strong law of large numbers for random elements in quasi-metric spaces. In *Trans. Seventh Prague Conference on Information Theory, Statistical Decision Functions, Random Processes and of the Eighth European Meeting of Statisticians* **A** 591–602. Czechoslovak Academy of Sciences, Prague. MR501230

DEPARTMENT OF STATISTICS
UNIVERSITY OF LEEDS
LEEDS LS2 9JT
UNITED KINGDOM
E-MAIL: k.v.mardia@leeds.ac.uk

DEPARTMENT OF MATHEMATICS
AND STATISTICS
TEXAS TECH UNIVERSITY
LUBBOCK, TEXAS 79409-1042
USA
E-MAIL: victor.patrangenaru@ttu.edu